\documentclass[11pt]{article}

\usepackage[T1]{fontenc}
\usepackage{lmodern}
\usepackage{microtype}
\usepackage{amsmath,amssymb,amsthm,mathtools}
\usepackage{bm}
\usepackage{enumitem}
\usepackage{mathrsfs}

\usepackage[margin=1.15in]{geometry}
\usepackage{graphicx}

\usepackage{hyperref}
\usepackage[nameinlink,capitalise,noabbrev]{cleveref}

\hypersetup{
    colorlinks=true,
    linkcolor=blue,
    citecolor=blue,
    urlcolor=blue
}

\newtheorem{theorem}{Theorem}[section]
\newtheorem{lemma}[theorem]{Lemma}

\theoremstyle{definition}
\newtheorem{definition}[theorem]{Definition}

\theoremstyle{remark}
\newtheorem{remark}[theorem]{Remark}

\title{Solvers for a Class of Coupled Forward-Backward Systems Under General Filtrations}
\author{
    Matthew Keyworth\\
    Department of Statistics, University of Warwick\\
    \texttt{matthew.keyworth@warwick.ac.uk}
}
\date{}

\begin{document}

\maketitle

\begin{abstract}
We consider a class of coupled-forward backward systems defined on a general probability space assumed only to satisfy the usual conditions. Under our assumptions we derive two solvers: a martingale solver (inspired by the deep FBSDE solver studied in recent literature), and a min-max solver leveraging a Hilbert space orthogonality condition to enforce the martingale property. We obtain a posteriori error bounds for both solvers. The min-max solver avoids the need to approximate the class of square integrable martingales directly (which in practice is usually done by parameterizing a stochastic integral representation), at the expense of weakening the topology with respect to which one should expect convergence of an approximate solution. 
\end{abstract}
\section{Introduction}
In this paper we study solution methods for a class of coupled forward-backward systems defined on a general probability space assumed only to satisfy the usual conditions. Existence and uniqueness results for the class of coupled forward-backward systems studied in this paper have previously been obtained in \cite{KeyworthPaper:2026}. In the case of coupled FBSDEs (typically with Brownian driver) methods of numerical solution have been studied in the literature. A modern approach that has been successful even in high dimensional problems is the deep FBSDE method (see \cite{PengPaper:2020}, \cite{PengPaper:2022},\cite{WangPaper:2022}, \cite{EvangelosPaper:2020}, \cite{EvangelosPaper:2021}, \cite{DaiPaper:2021}, \cite{LiuPaper:2021}) which evolved from the deep BSDE method for the numerical solution of BSDEs (see \cite{HanPaper:2018} \cite{WeinanPaper:2017}). In the deep FBSDE method, the solution of the coupled FBSDE is obtained as the unique solution of a variational problem involving the minimization of only a terminal error condition subject to forward dynamics that in practice are typically parameterized using deep neural networks. Minimizing the terminal error condition forces the parameterized forward dynamics to approximate the solution of the coupled FBSDE. In the literature some work has been done to examine the convergence of deep FBSDE method for numerical solution of certain coupled FBSDEs with Brownian driver (see \cite{HanPaper:2020}, \cite{JiangPaper:2021}, \cite{AnderssonPaper:2022}) but these all assume the existence of a so called decoupling field (see \cite{ProtterPaper:1994}, \cite{PardouxPaper:1999}, \cite{ZhangBook:2017}) which we cannot assume. In \cite{KeyworthPaper:2026} solving the coupled-forward backward system was reformulated equivalently as a fixed point problem for a certain nonlinear operator on a suitable Hilbert space. In this paper this reformulation again proves useful for studying the properties of two solvers that can then be discretized and implemented as stable numerical solvers. By stability, we mean that the solver objective value controls the distance between the approximate and true backward solution to the coupled forward-backward system in a certain norm. The first solver, the martingale solver, involves reformulating the coupled forward-backward system as a variational problem over the space of square integrable martingales (inspired by the deep FBSDE method recently studied in the literature). The second solver, the Hilbert space orthogonal min-max solver, involves reformulating the problem as a min-max game on a suitable Hilbert space. The numerical implementation of the martingale based solver would require the class of square integrable martingales to be parameterized in some way. For deep FBSDE solvers previously studied in the literature where the backward equation is a BSDE with a Brownian driver, the parameterization is achieved by approximating the initial condition and integrand of the stochastic integral using neural networks. This method can be applied in the case of our coupled forward-backward system if the filtration possess the martingale representation property with respect to a collection of known martingales that can be sampled. Our Hilbert space orthogonal min-max solver does not require any such restriction and is suitable for general filtrations. This comes at the expense of weakening the topology with respect to which one should expect convergence of the approximate solution to the true solution. We examine discretized versions of our solvers and prove suitable error bounds. In particular we show that the error of the approximate solution to the true solution in a certain norm is bounded by the discretized objective function value, subject to error terms coming from approximation error (that can be controlled by the choice of approximating class of functions), and time discretization error. In \cref{sec:preliminaries} we recall some basic facts and establish our notation. In \cref{sec:problem_formulation} we define our coupled forward-backward system and state our technical assumptions, and in \cref{sec:reformulation_as_fp} we recap the useful fixed point problem reformulation presented in \cite{KeyworthPaper:2026}. In \cref{sec:main_results} we detail our main results and how they fit together. In \cref{sec:proofs_main_results} we work through proofs of our main results.

\section{Preliminaries and Notation} \label{sec:preliminaries}
Define the following sets of stochastic processes:
\begin{equation}
    \begin{aligned}
         \mathcal{L}([0,T]; \mathbb{R}^d) &:=  \left\{\text{$X:\Omega\times[0,T]\rightarrow \mathbb{R}^d: \ X$  is $\mathcal{F}_T\otimes\mathcal{B}([0,T])$-measurable}\right\} \\
         \mathcal{L}_{\mathbb{F}}([0,T]; \mathbb{R}^d) &:= \left\{\text{$X \in \mathcal{L}([0,T]; \mathbb{R}^d)$ and $\mathbb{F}$-adapted}\right\} \\
        \mathcal{L}_{\mathbb{F}}^D([0,T]; \mathbb{R}^d) &:= \left\{X \in  \mathcal{L}_{\mathbb{F}}([0,T]; \mathbb{R}^d) : \text{$t\rightarrow X_t$ c\`adl\`ag $\mathbb{P}\text{-a.s.}$}\right\}\\
        \mathcal{L}_{\mathbb{F}}^C([0,T]; \mathbb{R}^d) &:= \left\{X \in  \mathcal{L}_{\mathbb{F}}^D([0,T]; \mathbb{R}^d) : \text{$t\rightarrow X_t$ continuous $\mathbb{P}\text{-a.s.}$}\right\} \\
        \mathcal{M}^2([0,T]; \mathbb{R}^d) &:= \left\{X \in \mathcal{L}^{D}_{\mathbb{F}}([0,T]; \mathbb{R}^d): \text{ $X$ is an $\mathbb{F}$-martingale with } \mathbb{E}[||X_T||_2^2] < \infty\right\} 
    \end{aligned}
\end{equation}
Now for any subset $A\subseteq \mathcal{L}([0,T]; \mathbb{R}^d)$ and any pseudo-metric $d:A\times A\rightarrow \mathbb{R}_+$ defining the equivalence relation:
\[X\sim_d Y \iff d(X, Y) = 0, \quad X,Y \in A\]
define the equivalence classes:
\begin{align*}
    [X]_d := \{Y \in A: Y\sim_d X\}, \quad \forall X \in A
\end{align*}
The pseudo-metrics:
\begin{align*}
    d_p(X, Y) := \mathbb{E}\left[\int_0^T ||X_t- Y_t ||_p^p\,ds\right]^{1/p}, \quad d_S(X, Y) = \mathbb{E}\left[\sup_{t \in [0,T]}||X_t - Y_t||_2^2\right]^{1/2}
\end{align*}
induce $d\mathbb{P}\otimes dt\text{-a.e.}$ equivalence and equivalence up to indistinguishability respectively. We thus define the following sets of equivalence classes:
\begin{equation}\label{eq:process_spaces}
    \begin{aligned}
         \mathbb{L}^p([0,T]; \mathbb{R}^d) &:= \left\{ [X]_{d_p} : X \in \mathcal{L}([0,T]; \mathbb{R}^d), \ d_p(X, 0) < \infty\right\}, \quad p \in \mathbb{N} \\
        \mathbb{L}^p_{\mathbb{F}}([0,T]; \mathbb{R}^d) &:= \left\{[X]_{d_p} : X \in \mathcal{L}_{\mathbb{F}}([0,T]; \mathbb{R}^d), d_p(X, 0) <\infty\right\}, \quad p \in \mathbb{N} \\
        \mathbb{S}^2_{\mathbb{F}}([0,T]; \mathbb{R}^d) &:= \left\{[X]_{d_S}: X \in \mathcal{L}^D_{\mathbb{F}}([0,T]; \mathbb{R}^d), d_S(X, 0) < \infty\right\} \\
        \mathbb{S}^{2,c}_{\mathbb{F}}([0,T]; \mathbb{R}^d) &:= \left\{[X]_{d_S}: X \in \mathcal{L}^C_{\mathbb{F}}([0,T]; \mathbb{R}^d), d_S(X, 0) < \infty\right\} \\
         \mathbb{M}^2_{\mathbb{F}}([0,T]; \mathbb{R}^d) &:= \left\{[X]_{d_S}: X \in \mathcal{M}^2([0,T]; \mathbb{R}^d), \ d_S(X, 0) < \infty \right\} 
    \end{aligned}
\end{equation}
For a set of equivalence classes $\mathbb{A}$ with respect to the equivalence relation $\sim_d$, we will use the usual abuse of notation and take $X \in \mathbb{A}$ to mean a representer of the equivalence class $[X]_d \in \mathbb{A}$.
We recall some basic facts about some important spaces that will be used in the subsequent sections. Equipped with the inner product:
\begin{align*}
    \langle X, Y \rangle &:= \mathbb{E}\left[\int_0^T\langle X_t, Y_t\rangle\,dt\right]
\end{align*}
the spaces $\mathbb{L}^2([0,T]; \mathbb{R}^d)$ and $\mathbb{L}^2_{\mathbb{F}}([0,T]; \mathbb{R}^d)$ are Hilbert spaces. Equipped with the norm induced by $d_S$, the spaces $\mathbb{S}^{2,c}_{\mathbb{F}}([0,T]; \mathbb{R}^d)$ and $\mathbb{S}^{2}_{\mathbb{F}}([0,T]; \mathbb{R}^d)$ are Banach spaces. The space $\mathbb{M}^{2}_{\mathbb{F}}([0,T]; \mathbb{R}^d)$ is Hilbert when equipped with the inner product:
\begin{align*}
    \langle M, M'\rangle = \mathbb{E}\left[M_T\cdot M'_T\right]
\end{align*}
and is a Banach subspace of $\mathbb{S}^{2}_{\mathbb{F}}([0,T]; \mathbb{R}^d)$ by norm equivalence:
 \begin{align*} 
    \mathbb{E}\left[||M_T||_2^2\right]^{1/2} \leq \mathbb{E}\left[\sup_{t \in [0,T]}||M_t||_2^2\right]^{1/2} \leq 2\mathbb{E}\left[||M_T||_2^2\right]^{1/2}
\end{align*}
which follows from Doob's inequality. By virtue of the usual conditions each $\mathbb{F}$-adapted martingale admits a c\`adl\`ag modification. We make use of the following canonical embeddings: 
\begin{equation}
    \begin{aligned}
        \mathbb{S}^2_{\mathbb{F}}([0,T]; \mathbb{R}^d)  \hookrightarrow\mathbb{L}^2_{\mathbb{F}}([0,T]; \mathbb{R}^d)  \xhookrightarrow{iso} \mathbb{L}^2([0,T]; \mathbb{R}^d) \\
    \end{aligned}
\end{equation}
Taking $i_1$ to be the canonical continuous embedding $ \mathbb{S}^2_{\mathbb{F}}([0,T]; \mathbb{R}^d)  \xhookrightarrow{i_1}\mathbb{L}^2_{\mathbb{F}}([0,T]; \mathbb{R}^d) $ defined by 
\begin{align*}
    i_1:\mathbb{S}^{2}_{\mathbb{F}}([0,T]; \mathbb{R}^d)\rightarrow\mathbb{L}^2_{\mathbb{F}}([0,T]; \mathbb{R}^d), \quad i_1([X]_{d_S}) = [X]_{d_2}
\end{align*}
we denote by $\mathcal{A}$ its image $\mathcal{A} := i_1(\mathbb{S}^{2}_{\mathbb{F}})\subseteq \mathbb{L}^2_{\mathbb{F}}$. Thus:
\begin{align} \label{eq:def_mathcal_A}
    \mathcal{A} &= \{[X]_{d_2} \in \mathbb{L}^2_{\mathbb{F}}([0,T]; \mathbb{R}^d): \text{$[X]_{d_2}$ admits an $\mathbb{F}$-adapted representative that is c\`adl\`ag $\mathbb{P}\text{-a.s.}$}\}
\end{align}
Taking $i_2$ to be the canonical isometric embedding $\mathbb{L}^2_{\mathbb{F}}([0,T]; \mathbb{R}^d)\xhookrightarrow{i_2}\mathbb{L}_{2}([0,T]; \mathbb{R}^d)$, defined by:
\begin{align*}
    i_2:\mathbb{L}^{2}_{\mathbb{F}}([0,T]; \mathbb{R}^d)\rightarrow\mathbb{L}^2([0,T]; \mathbb{R}^d), \quad i_2([X]_{d_2}) = [X]_{d_2}
\end{align*}
then we have:
\begin{align*}
       i_2(\mathbb{L}^2_{\mathbb{F}}) &=\left\{[X]_{d_2} \in \mathbb{L}^2([0,T]; \mathbb{R}^d) : \text{$[X]_{d_2} $ admits an $\mathbb{F}$-adapted representative}\right\} 
\end{align*}
Since $\mathbb{L}^2_{\mathbb{F}}$ is a Hilbert space and $i_2$ an isometry, $i_2(\mathbb{L}^2_{\mathbb{F}})\subseteq \mathbb{L}^2$ is a Hilbert subspace with $\mathbb{L}^2_{\mathbb{F}}\cong i_2(\mathbb{L}^2_{\mathbb{F}})$, and we have the orthogonal direct sum decomposition $\mathbb{L}^2 = i_2(\mathbb{L}^2_{\mathbb{F}})\oplus i_2(\mathbb{L}^2_{\mathbb{F}})^\perp$ and we define $\pi:\mathbb{L}^2\rightarrow i_2(\mathbb{L}^2_{\mathbb{F}})$ to be the orthogonal projection operator. Since $\mathcal{A},\ i_2(\mathbb{L}^2_{\mathbb{F}})$ are sets of equivalence classes, from now on we implicitly fix a representative. For $X \in \mathcal{A}$ we select the unique (up to indistinguishability) $\mathbb{F}$-adapted $\mathbb{P}\text{-a.s.}$ c\`adl\`ag representative. For and for $X \in i_2(\mathbb{L}^2_{\mathbb{F}})$ we select an $\mathbb{F}$-adapted representative.

\section{Problem Formulation and Assumptions} \label{sec:problem_formulation}
With $T \in (0,\infty)$ let $(\Omega, \mathcal{F}, \mathbb{F}:=(\mathcal{F}_t)_{t \in [0,T]}, \mathbb{P})$ be a filtered probability space satisfying the usual conditions. Let $g:\Omega \times \mathbb{R}^d\rightarrow \mathbb{R}^d$, $b:\Omega\times [0,T] \times \mathbb{R}^d\rightarrow \mathbb{R}^d$ and $\bar{\rho}:\Omega\times [0,T]\times \mathbb{R}^d\rightarrow U$ for convex $U\subseteq \mathbb{R}^d$. As is the usual convention of notation we mostly suppress the dependence of functions $g, b, \rho$ on $\omega$. We study solution methods for the following coupled forward-backward system:
\begin{equation} \label{eq:fbs_general}
    \begin{aligned}
        v_t &= v_0 + \int_0^t\bar{\rho}(s,\lambda_s)\,ds, \quad t \in [0,T], \quad v_0 \in \mathbb{R}^d\\
        \lambda_t &= \mathbb{E}\left[g(v_T) + \int_t^Tb(s, v_s)\,ds \ \bigg| \mathcal{F}_t \right], \quad \mathbb{P}\text{-a.s.}, \quad t \in [0,T] 
    \end{aligned}
\end{equation}
\begin{definition}[Solution Concept]\label{def:solution_concept}
We say that any pair $(v, \lambda) \in \mathbb{S}_{\mathbb{F}}^{2,c}([0,T]; \mathbb{R}^d)\times \mathbb{S}^2_{\mathbb{F}}([0,T]; \mathbb{R}^d)$ satisfying \cref{eq:fbs_general} is a solution.  We say that the forward-backward system has a unique solution iff for any two solutions $(v,\lambda)$ and $(v',\lambda')$ we have that:
\begin{align*}
    \mathbb{P}(v_t = v'_t; 0\leq t\leq T) = 1, \quad \mathbb{P}(\lambda_t = \lambda'_t; 0\leq t\leq T)= 1
\end{align*}
\end{definition}
We examine solution methods for \cref{eq:fbs_general} under the following assumptions:
\paragraph{Integrability Assumptions.}
With $v^\rho_t = v_0 + \int_0^t \rho_s\,ds, \ t \in [0,T]$:
\begin{enumerate}[label=(A\arabic*),  ref=A\arabic*]
 \item \label{ass:A1} $g(v_T^\rho) \in L^2(\Omega, \mathcal{F}_T, \mathbb{P}; \mathbb{R}^d), \ \forall \rho \in \mathbb{L}^2_{\mathbb{F}}([0,T]; \mathbb{R}^d)$
    \item \label{ass:A2} $b^\rho \in \mathbb{L}^2_{\mathbb{F}}([0,T]; \mathbb{R}^d), \forall \rho \in \mathbb{L}^2_{\mathbb{F}}([0,T]; \mathbb{R}^d)$ with $b^\rho_t := b(t,v^\rho_t), \ t \in [0,T]$
\item \label{ass:I1} $\sup_{t \in [0,T]}\mathbb{E}[||b(t, 0)||_2^2]  < \infty, \quad \sup_{t \in [0,T]}\mathbb{E}[||\bar{\rho}(t, 0)||_2^2]  < \infty$
\end{enumerate}

\paragraph{Monotonicity Assumptions.}
\begin{enumerate}[label=(M\arabic*),  ref=M\arabic*]
    \item \label{ass:B1} $\langle g(\omega, v) - g(\omega, v'), v - v'\rangle \geq 0, \quad \forall v, v' \in \mathbb{R}^d, \quad \mathbb{P}\text{-a.s.}$
    \item \label{ass:B2} $\exists \xi \in \mathbb{L}^1_{\mathbb{F}}([0,T]; \mathbb{R}), \ \xi_t > 0, \ \mathbb{P}\text{-a.s.}, \text{ a.e. } t \in [0,T]$ such that: 
    \[\langle b(t,\omega, v)-b(t,\omega, v'), v-v'\rangle \geq \xi(t,\omega)||v - v'||_2^2, \quad \forall v, v' \in \mathbb{R}^d, \ \mathbb{P}\text{-a.s.}, \text{ a.e. } t \in [0,T]\]
    \item \label{ass:B3}  $\exists\mu > 0$ such that $\forall \lambda, \lambda' \in \mathbb{R}^d:$
    \begin{align*}
        &\langle\bar{\rho}(t,\omega, \lambda) - \bar{\rho}(t,\omega, \lambda'), \lambda - \lambda'\rangle  \leq -\mu||\bar{\rho}(t,\omega, \lambda) - \bar{\rho}(t,\omega, \lambda')||_2^2, \quad \mathbb{P}\text{-a.s.}, \ \text{a.e.} \ t \in [0,T]  
    \end{align*}
\end{enumerate}
\paragraph{Global Lipschitz Assumptions.}
\begin{enumerate}[label=(A\arabic*), start=3, ref=A\arabic*]
    \item \label{ass:A5} There exist constants $L_g, L_b > 0$ such that:
    \begin{align*}
        ||g(\omega, v) - g(\omega, v')||_2 &\leq L_g||v - v'||_2, \quad v,v' \in \mathbb{R}^d, \quad \mathbb{P}\text{-a.s.} \\
        ||b(t,\omega, v) - b(t,\omega, v')||_2 &\leq L_b||v - v'||_2, \quad v,v' \in \mathbb{R}^d, \quad \mathbb{P}\text{-a.s.}, \ \text{a.e. } t \in [0,T]
    \end{align*}
\end{enumerate}
Note that by \ref{ass:B3} with $L_\rho = 1/\mu$ we obtain the global Lipschitz condition:
\begin{align*}
    ||\rho(t,\omega, \lambda) - \rho(t,\omega, \lambda')||_2 &\leq L_\rho||\lambda - \lambda'||_2, \quad \lambda, \lambda' \in \mathbb{R}^d, \quad \mathbb{P}\text{-a.s.}, \ \text{a.e. } t \in [0,T] 
\end{align*}

\section{Reformulation as a Fixed Point Problem}\label{sec:reformulation_as_fp}
Before summarizing the main results of this paper, we first recap some results we will require that were previously obtained in \cite{KeyworthPaper:2026}. In particular, defining the set:
\begin{align*}
    \mathcal{U}[0,T]:=\{\rho \in \mathcal{L}_{\mathbb{F}}([0,T]; \mathbb{R}^d): \ \rho(t,\omega) \in U, \ \mathbb{P}\text{-a.s.} \ \text{ a.e } t \in [0,T]\}
\end{align*}
which modulo $d\mathbb{P}\otimes dt\text{-a.e.}$ equivalence becomes the convex subset $\mathbb{U}[0,T]\subseteq \mathbb{L}^2_{\mathbb{F}}([0,T]; \mathbb{R}^d)$, examine the set of decoupled forward-backward systems:
\begin{equation} \label{eq:weakly_coupled_fbs}
    \begin{aligned}
        v_t &= v_0 + \int_0^t\rho_s\,ds, \quad \rho \in \mathbb{U}[0,T] \\
        \lambda_t^\rho &= \mathbb{E}\left[g(v_T) + \int_t^Tb(s,v_s)\,ds \ \bigg|\mathcal{F}_t\right], \quad \mathbb{P}\text{-a.s.}, \quad t \in [0,T]    
    \end{aligned}
\end{equation}
and define the strictly monotone operator $A:\mathbb{L}^2_{\mathbb{F}}([0,T]; \mathbb{R}^d)\rightarrow \mathbb{L}^2_{\mathbb{F}}([0,T]; \mathbb{R}^d)^*$:
\begin{equation}\label{eq:strictly_monotone_operator}
    \begin{aligned}
        A(\rho)(\varphi) &:= \mathbb{E}\left[g(v_T)\cdot v_T' + \int_0^Tb(t,v_t)\cdot v_t'\,dt\right] \\
        v_t &= v_0 + \int_0^t\rho_s\,ds, \quad  t \in [0,T], \quad \rho \in \mathbb{L}^2_{\mathbb{F}}([0,T]; \mathbb{R}^d)\\
        v_t' &= \int_0^t\varphi_s\,ds, \quad t \in [0,T], \quad \varphi \in \mathbb{L}^2_{\mathbb{F}}([0,T]; \mathbb{R}^d)
    \end{aligned}
\end{equation}
Then for each $\forall \rho \in \mathbb{U}[0,T]$, $\lambda^\rho \in \mathcal{A}$ uniquely satisfies the following duality relation:
\begin{equation} \label{eq:duality_rel}
    \begin{aligned}
        A(\rho)(\varphi) = \langle \lambda^\rho, \varphi\rangle_{\mathbb{L}^2_{\mathbb{F}}}, \quad \forall \varphi \in \mathbb{L}^2_{\mathbb{F}}([0,T]; \mathbb{R}^d), \ \forall \rho \in \mathbb{U}[0,T]
    \end{aligned}
\end{equation}
We can now re-frame the coupled forward-backward system of \cref{eq:fbs_general} as a fixed point problem on the Hilbert space $\mathbb{L}^2_{\mathbb{F}}([0,T]; \mathbb{R}^d)$. In particular, defining a Nemytskii map $\mathcal{N}:\mathbb{L}^2_{\mathbb{F}}([0,T]; \mathbb{R}^d)\rightarrow \mathbb{U}[0,T]$ such that:
\begin{equation}\label{eq:nemytskii_map}
    \begin{aligned}
        \mathcal{N}(\lambda)(t,\omega) = \bar{\rho}(t,\omega, \lambda(t,\omega)), \quad (t,\omega) \in [0,T]\times \Omega, \quad \lambda \in \mathbb{L}^2_{\mathbb{F}}([0,T]; \mathbb{R}^d)    
    \end{aligned}
\end{equation}
then from \ref{eq:weakly_coupled_fbs} we have the associated set of decoupled forward-backward systems:
\begin{equation} \label{eq:decoupled_fbs_parameterised}
    \begin{aligned}
        v_t &= v_0 + \int_0^t\mathcal{N}(\lambda)_s\,ds, \quad t \in [0,T], \quad \lambda \in \mathbb{L}^2_{\mathbb{F}}([0,T]; \mathbb{R}^d) \\
        \lambda_t^{\mathcal{N}(\lambda)} &= \mathbb{E}\left[g(v_T)+ \int_t^Tb(s,v_s)\,ds \ \bigg| \mathcal{F}_t\right], \quad \mathbb{P}\text{-a.s.}, \quad t \in [0,T] 
    \end{aligned}
\end{equation}
Using the duality relation of \ref{eq:duality_rel}, we have that  $\lambda^{\mathcal{N}(\lambda)} \in \mathcal{A}$ uniquely satisfies:
\begin{align}\label{eq:duality_rel_coupled}
    A(\mathcal{N}(\lambda))(\varphi) = \langle  \lambda^{\mathcal{N}(\lambda)}, \varphi \rangle_{\mathbb{L}^2_{\mathbb{F}}}, \quad \forall \varphi \in \mathbb{L}^2_{\mathbb{F}}([0,T]; \mathbb{R}^d), \ \lambda \in \mathbb{L}^2_{\mathbb{F}}([0,T]; \mathbb{R}^d)
\end{align}
Denoting the Riesz map by $R:\mathbb{L}^2_{\mathbb{F}}([0,T]; \mathbb{R}^d)\rightarrow \mathbb{L}^2_{\mathbb{F}}([0,T]; \mathbb{R}^d)^*$, we can define an important self map $\Phi:\mathcal{A}\rightarrow \mathcal{A}$ such that:
\begin{align} \label{eq:non_lin_phi_map}
    \lambda^{\mathcal{N}(\lambda)} &= \Phi(\lambda) := (R^{-1}\circ A \circ \mathcal{N})(\lambda), \quad \lambda \in \mathcal{A}
\end{align}
Defining the set of fixed points $\operatorname{Fix}\Phi := \{\lambda \in \mathcal{A} : \lambda = \Phi(\lambda)\}$ we have that:
\begin{align*}
    (v^\lambda,\lambda) \in \mathbb{S}^{2,c}_{\mathbb{F}}([0,T]; \mathbb{R}^d)\times \mathbb{S}^{2}_{\mathbb{F}}([0,T]; \mathbb{R}^d) \text{ solves \cref{eq:fbs_general}} \iff \lambda \in \operatorname{Fix}\Phi
\end{align*}
Therefore finding a solution to the coupled forward-backward system of \cref{eq:fbs_general} can be done by solving the fixed point problem to obtain $\lambda^* \in \operatorname{Fix}\Phi$ and then feeding $\lambda^*$ into the forward equation to obtain the solution pair $(v^*, \lambda^*)$. This re-framing was used in \cite{KeyworthPaper:2026} to prove existence and uniqueness results for \cref{eq:fbs_general}.
\section{Main Results} \label{sec:main_results}
Here we relay the main results of the paper and how they fit together. 
\subsection{Fixed Point Stability and Solvers}
In this paper we use the conceptual re-framing discussed in \cref{sec:reformulation_as_fp} to develop two solvers for \cref{eq:fbs_general}. In particular let $(v^*, \lambda^*) \in \mathbb{S}^{2,c}_{\mathbb{F}}([0,T]; \mathbb{R}^d)\times \mathbb{S}^{2}_{\mathbb{F}}([0,T]; \mathbb{R}^d)$ solve \cref{eq:fbs_general}, then we have $\lambda^* \in \operatorname{Fix}\Phi$ where $\operatorname{Fix}\Phi\subseteq \mathcal{A}$ is a singleton. Consequently, with respect to the norm $||\cdot||$ induced by $d_2$ and $d_S$ we have that $\lambda^*$ uniquely solves the fixed point problem:
\begin{align*}
    0 = ||(I-\Phi)(\lambda^*)|| = \inf_{\lambda \in \mathcal{A}}||(I-\Phi)(\lambda)|| 
\end{align*}
One would hope therefore to solve the coupled system of \cref{eq:fbs_general} by minimizing the norm of the fixed point residual $r(\lambda):=(I-\Phi)(\lambda)$ over the set $\mathcal{A}$, which would force the decoupled system of \cref{eq:decoupled_fbs_parameterised} to become coupled. However, given a sequence $(\lambda_n)_{n \in \mathbb{N}} \subseteq \mathcal{A}$ with $||r(\lambda_n)||\rightarrow 0$, we cannot in general claim that this implies convergence $\lambda_n\rightarrow\lambda^*$. 
Crucially we can only claim that such a solver would be useful numerically if we can prove the additional property that.
\begin{align*}
    ||(I-\Phi)(\lambda^n)||\rightarrow 0 \ \implies ||\lambda^n - \lambda^*|| \rightarrow 0
\end{align*}
with respect to norms $||\cdot||$ induced by $d_2$ and/ or $d_S$. This leads us to a natural definition of solver stability. This leads us to our first main result.
\begin{theorem}[Fixed Point Stability]\label{lemma:fixed_point_stability}
    Let $\lambda^* \in \operatorname{Fix}\Phi$, then there exists $C > 0$ such that:
    \begin{enumerate}
        \item $  d_2(\lambda, \lambda^*) \leq Cd_2(\Phi(\lambda), \lambda), \quad \forall \lambda \in \mathcal{A}$
        \item $d_S(\lambda, \lambda^*) \leq Cd_S(\Phi(\lambda), \lambda), \quad \forall \lambda \in \mathcal{A}$ 
    \end{enumerate}
\end{theorem}
The above fixed point stability result helps us to investigate the properties of two theoretical solvers: a martingale based solver (which takes inspiration from the deep FBSDE method studied in recent literature), and a solver formulated as a min-max game.

\subsection{Solvers}
We will claim that any function that dominates the fixed point residual in a certain norm with a unique minimiser at $\lambda^* \in \operatorname{Fix}\Phi$ is a theoretical solver. We will claim such a solver is stable if the function value controls the distance between the approximate solution and the true solution in a certain norm.
\begin{definition}[Solver / Solver Stability] \label{def:solver}
    We will call a function $J:\mathcal{A}\rightarrow \mathbb{R}_+$ a solver for the fixed point problem (and hence \cref{eq:fbs_general}) with respect to the norm $||\cdot||$ if $\lambda^* \in \operatorname{Fix}\Phi$ is the unique minimizer:
    \begin{align*}
        \inf_{\lambda \in \mathcal{A}}J(\lambda) = J(\lambda^*) = 0
    \end{align*}
    and there exist constants $C > 0$ and $\alpha > 0$ such that:
    \begin{align*}
        ||(I-\Phi)(\lambda)|| \leq  CJ(\lambda)^\alpha, \quad \forall \lambda \in \mathcal{A}
    \end{align*}
    We will call a solver $J:\mathcal{A}\rightarrow \mathbb{R}_+$ stable with respect to the norm $||\cdot||$ if there exist constants $C > 0$, $\alpha > 0$ such that:
    \begin{align*}
        ||\lambda^* - \lambda|| \leq CJ(\lambda)^\alpha
    \end{align*}
\end{definition}

\subsubsection{Martingale Solver}
Inspired by the Deep FBSDE solver defined in recent literature, we define the following variational problem on the space of square integrable martingales $\mathbb{M}^{2}_{\mathbb{F}}([0,T]; \mathbb{R}^d)$ subject to forward dynamics:
\begin{equation} \label{eq:martingale_solver}
    \begin{aligned}
        J(M) &= \mathbb{E}[||g(v_T) - \lambda_T||_2^2], \quad M \in \mathbb{M}^{2}_{\mathbb{F}}([0,T]; \mathbb{R}^d) \\
    &\text{s.t.} \\
    \lambda_t &= - \int_0^tb(s,v_s)\,ds + M_t \\
    v_s &= v_0 + \int_0^t\bar{\rho}(s,\lambda_s)\,ds 
    \end{aligned}
\end{equation}
\begin{theorem}[Martingale Solver]\label{theorem:martingale_solver_stable}
    The variational problem of \cref{eq:martingale_solver} 
    \begin{enumerate}
    \item is a solver in the sense of the \cref{def:solver} for the coupled forward backward system of \cref{eq:fbs_general} with respect to the norms induced by $d_S$ and $d_2$.  I.e. there exist $C_1, C_2 > 0$ such that:
    \begin{align*}
        d_S(\Phi(\lambda), \lambda) \leq C_2J(M)^{1/2}, \quad d_2(\Phi(\lambda), \lambda) \leq C_2J(M)^{1/2}, \quad M \in \mathbb{M}^{2}_{\mathbb{F}}([0,T]; \mathbb{R}^d)
    \end{align*}
    and there exists a unique $M^* \in \mathbb{M}^{2}_{\mathbb{F}}([0,T];\mathbb{R}^d)$ such that:
    \begin{align*}
        \inf_{M \in\mathbb{M}^{2}_{\mathbb{F}}}J(M) = J(M^*) = 0
    \end{align*}
    and where:
    \begin{align*}
        M^*_t &= \mathbb{E}\left[g(v^*_T) + \int_0^Tb(s,v^*_s)\,ds \ \bigg|\mathcal{F}_t\right], \quad \mathbb{P}\text{-a.s.}, \quad t \in [0,T] \\
        \lambda^*_t &= - \int_0^tb(s,v^*_s)\,ds + M^*_t, \quad t \in [0,T] \\
        v^*_t &= v_0 + \int_0^t\bar{\rho}(s,\lambda^*_s)\,ds, \quad t \in [0,T]
    \end{align*}
    and where the pair $(v^*, \lambda^*)$ solves \cref{eq:fbs_general}.
    \item the variational problem of \cref{eq:martingale_solver} is a stable in the sense of \cref{def:solver} with respect to norms induced by $d_S$ and $d_2$. I.e. there exists $C_1, C_2 > 0$ such that:
    \begin{align*}
         d_S(\lambda, \lambda^*) \leq C_1J(M)^{1/2}, \quad d_2(\lambda, \lambda^*) \leq C_2J(M)^{1/2}, \quad \forall M \in \mathbb{M}^{2}_{\mathbb{F}}([0,T]; \mathbb{R}^d)
    \end{align*}    
    \end{enumerate}
\end{theorem}

\subsubsection{Hilbert Space Orthogonal Min-Max Solvers}
We then define a solver that does not require the solution of a variational problem on $\mathbb{M}^{2}_{\mathbb{F}}([0,T]; \mathbb{R}^d)$ directly, as is the case for the martingale solver. Instead we can enforce the martingale property using a Hilbert space orthogonality condition. This naturally leads to a two player min-max game for which the unique zero is achieved when \cref{eq:fbs_general} is solved. Letting $J:\mathcal{A}\rightarrow \mathbb{R}$ such that:
\begin{equation} \label{eq:controlled_dynamics}
    \begin{aligned}
        J(\lambda) &= \mathbb{E}\left[||g(v_T) - \lambda_T||_2^2\right], \quad \lambda \in \mathcal{A} \\
        &\text{s.t.} \\
        v_t &= v_0 + \int_0^t\bar{\rho}(s,\lambda_s)\,ds, \quad t \in [0,T] \\
        M_t &= \lambda_t + \int_0^tb(s, v_s)\,ds, \quad t \in [0,T]
    \end{aligned}
\end{equation}
then for each $\lambda \in \mathcal{A}$ importantly $M \in \mathcal{A}$ but not necessarily in $\mathbb{M}^{2}_{\mathbb{F}}([0,T]; \mathbb{R}^d)$. Defining what we call the $\textit{targets}$, such that:
\begin{align*}
    (Y_t)_{0\leq t\leq T} \equiv g(v_T) + \int_0^Tb(s,v_s)\,ds
\end{align*}
and recalling the orthogonal projection operator $\pi:\mathbb{L}^2\rightarrow i_2(\mathbb{L}^2_{\mathbb{F}})$ we can introduce a constraint function $F : \mathcal{A}\rightarrow i_2(\mathbb{L}^2_{\mathbb{F}})\subseteq \mathbb{L}^2$ such that:
\begin{align*}
    F(\lambda) := \pi(Y-M) = \pi Y - \pi M = \pi Y - M
\end{align*}
We will see that solving \cref{eq:fbs_general} equates to solving the variational problem of \ref{eq:controlled_dynamics} subject to $F(\lambda) = 0$ acting as a martingale constraint, and that this formulation leads naturally to a min-max game. We then prove that the min-max game is a stable solver for \cref{eq:fbs_general} in the sense of \cref{def:solver}.
\begin{theorem}\label{theorem:unique_sol_const_opt}
There exists a unique feasible solution $\lambda^* \in \mathcal{A}$ to the constrained optimisation problem:
\begin{equation} \label{eq:const_opt}
    \begin{aligned}
        \inf_{\lambda \in \mathcal{A}}J(\lambda) \quad  \text{ s.t. } \quad F(\lambda) = 0
    \end{aligned}
\end{equation}
and that the pair $(v^*, \lambda^*)$ solves \cref{eq:fbs_general}.
\end{theorem}
\begin{theorem}[Min-Max Solvers] \label{theorem:min_max_game}
    Define the Lagrangian $L:\mathcal{A}\times \mathbb{L}^2_{\mathbb{F}}\rightarrow \mathbb{R}$ such that:
    \begin{align}\label{eq:lagrangian}
         L(\lambda, \varphi) := J(\lambda) + \langle \varphi, Y - M\rangle_{\mathbb{L}^2}, \quad \lambda \in \mathcal{A}, \quad \varphi \in \mathbb{L}^2_{\mathbb{F}}
    \end{align}
    and the regularized Lagrangian $\tilde{L}:\mathcal{A}\times \mathbb{L}^2_{\mathbb{F}}\rightarrow \mathbb{R}$ such that for $r > 0$:
    \begin{align}\label{eq:lagrangian_reg}
         \tilde{L}(\lambda, \varphi) := J(\lambda) + \langle \varphi, Y - M\rangle_{\mathbb{L}^2} - \frac{r}{2}||\varphi||_{\mathbb{L}^2_{\mathbb{F}}}^2, \quad \lambda \in \mathcal{A}, \quad \varphi \in \mathbb{L}^2_{\mathbb{F}}
    \end{align}
    Let $\mathcal{D}\subseteq \mathbb{L}^2_{\mathbb{F}}$ be a dense linear subspace and for $c > 0$ define $\mathcal{T}_c = \{\varphi \in \mathcal{D}: ||\varphi||_{\mathbb{L}^2_{\mathbb{F}}}\leq c\}$. Then:
    \begin{enumerate}
        \item The Lagrangian of \ref{eq:lagrangian} and regularized Lagrangian of \ref{eq:lagrangian_reg} are solvers in the sense of \cref{def:solver} with respect to the norm induced by $d_2$. I.e. there exist $C_1, C_2 > 0$ such that:
        \begin{align*}
            d_2(\Phi(\lambda), \lambda) \leq C_1\sup_{\varphi \in \mathcal{T}_c}L(\lambda, \varphi), \quad d_2(\Phi(\lambda), \lambda) \leq C_2\sup_{\varphi \in \mathcal{D}}\tilde{L}(\lambda, \varphi), \quad \forall \lambda \in \mathcal{A}
        \end{align*}
        and $\lambda^* \in \operatorname{Fix}\Phi$ is the unique zero. 
        \begin{align*}
            0 = \sup_{\varphi \in \mathcal{T}_c}L(\lambda^*, \varphi)= \inf_{\lambda \in \mathcal{A}}\sup_{\varphi \in \mathcal{T}_c}\tilde{L}(\lambda, \varphi), \quad 0 = \sup_{\varphi \in \mathcal{D}}L(\lambda^*, \varphi)= \inf_{\lambda \in \mathcal{A}}\sup_{\varphi \in \mathcal{D}}\tilde{L}(\lambda, \varphi),
        \end{align*}
        \item The Lagrangian of \ref{eq:lagrangian} and regularized Lagrangian of \ref{eq:lagrangian_reg} are stable solvers in the sense of \cref{def:solver} with respect to the norm induced by $d_2$. I.e. there exist $C_1, C_2 > 0$ such that:
        \begin{align*}
            d_2(\lambda, \lambda^*) \leq C_1\sup_{\varphi \in \mathcal{T}_c}L(\lambda, \varphi), \quad d_2(\lambda, \lambda^*)^2 \leq C_2\sup_{\varphi \in \mathcal{D}}L(\lambda, \varphi), \quad \forall \lambda \in \mathcal{A}
        \end{align*}
    \end{enumerate}
\end{theorem}

\subsection{Discretized Solvers}
 We study discretized versions of the martingale and min-max solvers and our fixed point stability results again prove useful.  We will show that the distance of the true solution $\lambda^* \in \operatorname{Fix}\Phi$ to our approximate solution is bounded by the value of the discretized cost plus error terms corresponding to time discretization error and approximation error (error due the limited capacity of the class of approximating functions within which we search for a solution). For $N \in \mathbb{N}$ we take a partition of $[0,T]$ such that $0 = t_0 < t_1 < \dots < t_N = T$ with $\Delta_N = T/N$, and as in \cite{JiangPaper:2021} define $\tau:[0,T]\rightarrow \{t_0, \dots, t_N\}$ such that $\tau(t):= \left\lfloor t(\Delta_N )^{-1} \right\rfloor\Delta_N$. 
\subsubsection{Martingale Solver}
Let $\mathscr{A}_M\subseteq \mathbb{M}^{2}_{\mathbb{F}}([0,T]; \mathbb{R}^d)$ be an approximating class and define the forward dynamics for the discretized problem:
\begin{equation}\label{eq:mart_solver_dynamics_tau}
    \begin{aligned}
        \lambda_t^N &= -\int_0^tb(\tau(s), v_{\tau(s)}^N)\,ds + M_t^N, \quad t \in [0, T], \quad M^N \in \mathscr{A}_M \\
    v_t^N &= v_0 + \int_0^t\bar{\rho}(\tau(s), \lambda_{\tau(s)}^N)\,ds, \quad t \in [0,T]    
    \end{aligned}
\end{equation}
Thus we have that:
\begin{align*}
    \lambda_{t_{i}}^N &= -\sum_{k=0}^{i-1}b(t_k, v_{t_k}^N)\,\Delta + M_{t_i}^N, \quad i=1, \dots, N, \quad \lambda_{t_0}^N = M_0^N \\
    v_{t_i}^N &= v_0 + \sum_{k=0}^{i-1}\bar{\rho}(t_k, \lambda_{t_k}^N)\,\Delta, \quad i=1, \dots, N, \quad v_{t_0}^N = v_0
\end{align*}
and we have a well defined forward scheme for the grid points of the partition:
\begin{align*}
    \lambda_{t_{i}}^N &= \lambda_{t_{i-1}}^N -b(t_{i-1}, v_{t_{i-1}}^N)\,\Delta + M_{t_i}^N - M_{t_{i-1}}^N, \quad i=1, \dots, N, \quad \lambda_{t_0}^N = M_0^N \\
    v_{t_i}^N &= v_{t_{i-1}}^N + \bar{\rho}(t_{i-1}, \lambda_{t_{i-1}}^N)\,\Delta, \quad i=1, \dots, N, \quad v_{t_0}^N = v_0
\end{align*}
We define the objective function for the discretized martingale solver as follows:
\begin{equation}\label{eq:discrete_martingale_solver}
    \begin{aligned}
        J_N(M^N) &= \mathbb{E}\left[||g(v_{t_N}^N) - \lambda_{t_N}^N||_2^2\right], \quad M^N \in \mathscr{A}_M\\
        &\text{s.t.} \\
        \lambda_t^N &= -\int_0^tb(\tau(s), v_{\tau(s)}^N)\,ds + M_t^N \\
        v_t^N &= v_0 + \int_0^t\bar{\rho}(\tau(s), \lambda_{\tau(s)}^N)\,ds
    \end{aligned}
\end{equation}
We will prove that the discretized problem of \ref{eq:discrete_martingale_solver} is a solver for \cref{eq:fbs_general} subject to approximation error depending on the choice of $\mathscr{A}_M$ and discretization error depending on $N \in \mathbb{N}$.
\begin{theorem}[Discretized Martingale Solver] \label{theorem:martingale_solver_discrete}
    There exists $C > 0$ such that the discretized martingale solver of \cref{eq:discrete_martingale_solver} satisfies the following:
    \begin{align*}
        d_S(\lambda^N, \lambda^*)^2 \leq C\left(J_N(M^N) + \delta_N\right), \quad \forall M^N \in \mathscr{A}_M, \quad \forall N \in \mathbb{N} 
    \end{align*}
    where $\delta_N \geq 0$ depends on $N$, and where:
    \begin{align*}
        0\leq \inf_{M^N \in \mathscr{A}_M}J_N(M^N) \leq \sup_{M \in \mathbb{M}^{2}_{\mathbb{F}}}d_S(M, \mathscr{A}_M)^2
    \end{align*}
\end{theorem}
\subsubsection{Min-Max Solvers}
Let $\mathscr{A}_{\lambda} \subseteq \mathcal{A}$ be an approximating class, and for any dense subspace $\mathcal{D}\subseteq\mathbb{L}^2_{\mathbb{F}}([0,T]; \mathbb{R}^d)$ and any linear subspace $\tilde{\mathcal{D}}\subseteq \mathcal{D}$ with $c > 0$ define:
\[\tilde{\mathcal{T}}_c:=\{\varphi \in \tilde{\mathcal{D}} : ||\varphi||_{\mathbb{L}^2_{\mathbb{F}}} \leq c\}\]
  Given $\lambda^N \in \mathscr{A}_{\lambda}$ define the forward dynamics and targets:
\begin{equation}\label{eq:discretized_dynamics_min_max_solver}
    \begin{aligned}
        v_t^N &= v_0 + \int_0^t\bar{\rho}(\tau(s),\lambda_{\tau(s)}^N)\,ds, \quad t \in [0,T] \\
        M_t^N &= \lambda_t^N + \int_0^tb(\tau(s), v_{\tau(s)}^N)\,ds, \quad t \in [0,T] \\
        (Y_t^N)_{0\leq t\leq T} &\equiv g(v_{t_N}^N) + \int_0^Tb(\tau(s),v_{\tau(s)}^N)\,ds    
    \end{aligned}
\end{equation}
Thus we have:
\begin{align*}
    v_{t_i}^N &= v_0 + \sum_{k=0}^{i-1}\bar{\rho}(t_k,\lambda_{t_k}^N)\,\Delta, \quad i = 1, \dots, N, \quad  v_{t_0}^N = v_0\\
    M_{t_i}^N &= \lambda_{t_i}^N + \sum_{k=0}^{i-1}b(t_k, v_{t_k}^N)\,\Delta, \quad i= 1, \dots, N, \quad  M_{t_0}^N = \lambda_{t_0}^N\\
    (Y_{t}^N)_{0\leq t\leq T} &\equiv g(v_{t_N}^N) + \sum_{k=0}^{N-1}b(t_k,v_{t_k}^N)\,\Delta
\end{align*}
and we have enforced the following dynamics at the grid points of the partition:
\begin{align*}
    v_{t_i}^N  &= v_{t_{i-1}}^N + \bar{\rho}(t_{i-1},\lambda_{t_{i-1}}^N)\,\Delta, \quad i = 1, \dots, N, \quad  v_{t_0}^N = v_0\\
    \lambda_{t_i}^N &=  \lambda_{t_{i-1}}^N - b(t_{i-1}, v_{t_{i-1}}^N)\,\Delta + M_{t_i}^N - M_{t_{i-1}}^N, \quad i= 1, \dots, N, \quad  \lambda_{t_0}^N = M_{t_0}^N 
\end{align*}
We define the objective functions for the discrete min-max solvers as follows:
\begin{equation} \label{eq:discrete_minmax_solvers}
    \begin{aligned}
        L_N(\lambda^N, \varphi) &= \mathbb{E}\left[||g(v_{t_N}^N) - \lambda_{t_N}^N||_2^2\right]^{1/2} + \mathbb{E}\left[\int_0^T\langle \varphi_{s}, Y_{s}^N - M_{s}^N\rangle\,ds\right], \quad \varphi \in \tilde{\mathcal{T}}_c \\
        \tilde{L}_N(\lambda^N, \varphi) &= \mathbb{E}\left[||g(v_{t_N}^N) - \lambda_{t_N}^N||_2^2\right] + \mathbb{E}\left[\int_0^T\langle \varphi_{s}, Y_{s}^N - M_{s}^N\rangle\,ds\right] - \frac{r}{2}\mathbb{E}\left[\int_0^T||\varphi_{s}||_2^2\,ds\right], \quad \varphi \in \tilde{\mathcal{D}} \\
        &\text{s.t.} \\
        v_t^N &= v_0 + \int_0^t\bar{\rho}(\tau(s),\lambda_{\tau(s)}^N)\,ds, \quad t \in [0,T] \\
    M_t^N &= \lambda_t^N + \int_0^tb(\tau(s), v_{\tau(s)}^N)\,ds, \quad t \in [0,T] \\
    (Y_t^N)_{0\leq t\leq T} &\equiv g(v_{t_N}^N) + \int_0^Tb(\tau(s),v_{\tau(s)}^N)\,ds 
    \end{aligned}
\end{equation}
We will prove that the discretized min-max problems of \cref{eq:discrete_minmax_solvers} are solvers for \cref{eq:fbs_general} subject to approximation error depending on the choice of approximating class $\mathscr{A}_\lambda$, test function class $\tilde{\mathcal{D}}$, and discretization error depending on $N \in \mathbb{N}$.
\begin{theorem}[Discretized Min-Max Solvers]\label{theorem:discretized_min_max_solvers}
   Suppose that there exists $C_{\lambda} > 0$ such that $||\lambda ||_{\mathbb{S}^2_{\mathbb{F}}}< C_{\lambda}, \ \forall \lambda \in \mathscr{A}_{\lambda}$. Then
    \begin{enumerate}
        \item there exist $C_1, C_2 > 0$ such that:
        \begin{equation}\label{eq:discrete_min_max_lag_bound}
        \begin{aligned}
            d_2(\lambda^N, \lambda^*) &\leq C_1\left(\sup_{\varphi \in \tilde{\mathcal{T}}_c}L_N(\lambda^N, \tilde{\varphi})+ \sup_{\varphi \in \mathcal{D}}d_2(\varphi, \tilde{\mathcal{D}}) + \delta_N^{1/2}\right) \\
            \inf_{\lambda^N \in \mathscr{A}_{\lambda}}\sup_{\tilde{\varphi} \in \tilde{\mathcal{T}}_c}L_N(\lambda^N, \tilde{\varphi}) &\leq C_2\sup_{M \in \mathbb{M}^{2}_{\mathbb{F}}}\inf_{\lambda^N \in \mathscr{A}_{\lambda}}d_S(M, M^N(\lambda^N)) 
        \end{aligned}
        \end{equation}
        \item there exist $C_1, C_2 > 0$ such that:
        \begin{equation}\label{eq:discrete_min_max_reg_lag_bound}
        \begin{aligned}
             d_2(\lambda^N, \lambda^*)^2 &\leq C_1\left(\sup_{\tilde{\varphi} \in \tilde{\mathcal{D}}}\tilde{L}_N(\lambda^N, \tilde{\varphi})
            +\sup_{\varphi \in \mathcal{D}}d_2(\varphi, \tilde{\mathcal{D}}) +\delta_N^{1/2}\right)\\
            \inf_{\lambda^N \in \mathscr{A}_{\lambda}}\sup_{\tilde{\varphi} \in \tilde{\mathcal{D}}}\tilde{L}_N(\lambda^N, \tilde{\varphi}) &\leq C_2\sup_{M \in \mathbb{M}^{2}_{\mathbb{F}}}\inf_{\lambda^N \in \mathscr{A}_{\lambda}}d_S(M, M^N(\lambda^N))^2
        \end{aligned}
    \end{equation}
    \end{enumerate}
    with $\delta_N\geq 0$ depending on $N$:
\end{theorem}
\section{Proofs of Main Results}\label{sec:proofs_main_results}
\subsection{Fixed Point Stability} \label{subsec:fixed_point_stability}
Before proving \cref{lemma:fixed_point_stability} we will prove a number of useful lemmas.
\begin{lemma}\label{lemma:d2_bound_control}
    Let $\lambda^* \in \operatorname{Fix}\Phi$, then
    \begin{align*}
        ||\Phi(\lambda) - \lambda||_{\mathbb{L}^2_{\mathbb{F}}} \geq \mu||\mathcal{N}(\lambda)- \mathcal{N}(\lambda^*)||_{\mathbb{L}^2_{\mathbb{F}}}, \quad \forall \lambda \in \mathcal{A}
    \end{align*}
    \begin{proof}
        Let $\rho = \mathcal{N}(\lambda)$ and $\rho^* = \mathcal{N}(\lambda^*)$. From the duality relation of \ref{eq:duality_rel}, the definition of $\Phi$, and monotonicity of $A$ we have:
        \begin{align*}
            A(\rho)(\rho - \rho^*) - A(\rho^*)(\rho - \rho^*)= \langle \Phi(\lambda) - \lambda^*, \rho - \rho^*\rangle_{\mathbb{L}^2_{\mathbb{F}}} \geq 0, \quad \forall \lambda \in \mathcal{A}
        \end{align*}
        Thus:
        \begin{align*}
            \langle \Phi(\lambda) - \lambda^*, \rho - \rho^*\rangle_{\mathbb{L}^2_{\mathbb{F}}} = \langle \Phi(\lambda) - \lambda + \lambda - \lambda^*, \rho - \rho^*\rangle_{\mathbb{L}^2_{\mathbb{F}}} \geq 0, \quad \forall \lambda \in \mathcal{A}
        \end{align*}
        By \ref{ass:B3} we obtain:
        \begin{align*}
            \langle \Phi(\lambda) - \lambda , \rho - \rho^*\rangle_{\mathbb{L}^2_{\mathbb{F}}} &\geq -\langle \lambda - \lambda^*, \rho - \rho^*\rangle_{\mathbb{L}^2_{\mathbb{F}}} \\
            &= -\mathbb{E}\left[\int_0^T(\lambda_t - \lambda^*_t)\cdot(\bar{\rho}(t,\lambda_t) - \bar{\rho}(t,\lambda^*_t))\,dt\right] \\
            &\geq \mu||\rho - \rho^*||_{\mathbb{L}^2_{\mathbb{F}}}^2, \quad \forall \lambda \in \mathcal{A}
        \end{align*}
        and by Cauchy-Schwarz we obtain the result:
        \begin{align*}
            ||\Phi(\lambda) - \lambda||_{\mathbb{L}^2_{\mathbb{F}}}   \geq \mu||\rho - \rho^*||_{\mathbb{L}^2_{\mathbb{F}}},\ \quad \forall \lambda \in \mathcal{A}
        \end{align*} 
    \end{proof}
\end{lemma}
\begin{lemma}\label{lemma:bounding_mart}
    For every $\lambda^1,\lambda^2 \in \mathcal{A}$ the non-negative process:
    \begin{align*}
        M_t &= \mathbb{E}\left[||g(v_T^1) - g(v_T^2)||_2 + \int_0^T||b(s,v_s^1) - b(s,v^2_s)||_2\,ds \ \bigg|\mathcal{F}_t\right], \quad t \in [0,T] \\
        v^i_t &= v_0 + \int_0^t\bar{\rho}(s,\lambda_s^i)\,ds, \quad t \in [0,T], \quad i \in \{1,2\}
    \end{align*}
    is a square integrable martingale with:
    \begin{align*}
         \mathbb{E}\left[\sup_{0\leq t\leq T}|M_t|^2\right]\leq 4\mathbb{E}\left[\left(||g(v_T^1) - g(v^2_T)||_2 + \int_0^T||b(s,v_s^1) - b(s,v^2_s)||_2\,ds\right)^2\right]
    \end{align*}
    \begin{proof}
        We have for each $t \in [0,T]$ that:
        \begin{align*}
            \mathbb{E}\left[|M_t|^2\right] &\leq 
             \mathbb{E}\left[\bigg|||g(v_T^1) - g(v_T^2)||_2 + \int_0^T||b(s,v_s^1) - b(s,v^2_s)||_2\,ds \bigg|^2\right] \\
             &\leq 2\left(\mathbb{E}\left[||g(v_T^1) - g(v_T^2)||_2^2\right] + T\mathbb{E}\left[ \int_0^T||b(s,v_s^1) - b(s,v^2_s)||_2^2\,ds \right]\right) \\
             &\leq 4\left(||g(v^1_T)||_{L^2_{\mathcal{F}_T}} + ||g(v^2_T)||_{L^2_{\mathcal{F}_T}} + T||b^1||_{\mathbb{L}^2_{\mathbb{F}}} + T||b^2||_{\mathbb{L}^2_{\mathbb{F}}}\right) < \infty
        \end{align*}
        where the RHS is finite by \ref{ass:A1} and \ref{ass:A2}.  Hence $M$ is a non-negative square integrable martingale, and by Doob's submartingale inequality we obtain:
        \begin{align*}
            \mathbb{E}\left[\sup_{0\leq t\leq T}|M_t|^2\right]  
            \leq 4\mathbb{E}\left[\left(||g(v_T^1) - g(v^2_T)||_2 + \int_0^T||b(s,v_s^1) - b(s,v^2_s)||_2\,ds\right)^2\right] 
        \end{align*}
    \end{proof}
\end{lemma}
\begin{lemma} \label{lemma:fixed_point_stability_bkwd_dS}
    Let $\lambda^* \in \operatorname{Fix}\Phi$, then there exists $C > 0$ such that:
    \begin{align*}
        d_S(\Phi(\lambda), \lambda^*) \leq Cd_2(\Phi(\lambda), \lambda), \quad \forall \lambda \in \mathcal{A}
    \end{align*}
    \begin{proof}
        For each $\lambda \in \mathcal{A}$ we have that
        \begin{align*}
            \Phi_t(\lambda) - \lambda_t^* = \mathbb{E}\left[g(v_T) - g(v^*_T) + \int_t^T(b(s,v_s) - b(s,v^*_s))\,ds \ \bigg|\mathcal{F}_t\right], \quad \mathbb{P}\text{-a.s.}, \quad t \in [0,T]
        \end{align*}
        Now define the non-negative process $M = (M_t)_{t \in [0,T]}$ such that:
        \begin{align*}
            M_t = \mathbb{E}\left[||g(v_T) - g(v^*_T)||_2 + \int_0^T||b(s,v_s) - b(s,v^*_s)||_2\,ds \ \bigg|\mathcal{F}_t\right], \quad t \in [0,T]
        \end{align*}
        which is a square integrable martingale by \cref{lemma:bounding_mart} satisfying:
        \begin{align*}
         \mathbb{E}\left[\sup_{0\leq t\leq T}|M_t|^2\right]\leq 4\mathbb{E}\left[\left(||g(v_T) - g(v^*_T)||_2 + \int_0^T||b(s,v_s) - b(s,v^*_s)||_2\,ds\right)^2\right]
        \end{align*}
        Clearly we have that:
        \begin{align*}
            ||\Phi_t(\lambda) - \lambda^*_t||_{2}^2 \leq |M_t|^2, \quad \forall t \in [0,T]
        \end{align*}
        Thus we have that:
        \begin{align*}
            \mathbb{E}\left[\sup_{0\leq t\leq T}||\Phi(\lambda)_t - \lambda^*_t||_2^2\right] \leq 4\mathbb{E}\left[\left(||g(v_T) - g(v^*_T)||_2 + \int_0^T||b(s,v_s) - b(s,v^*_s)||_2\,ds\right)^2\right] \\
            \leq 8\mathbb{E}\left[||g(v_T) - g(v^*_T)||_2^2 + T\int_0^T||b(s,v_s) - b(s,v^*_s)||_2^2\,ds\right]
        \end{align*}
        and by \ref{ass:A5} we obtain:
        \begin{align*}
             \mathbb{E}\left[\sup_{0\leq t\leq T}||\Phi(\lambda)_t - \lambda^*_t||_2^2\right] &\leq 8\mathbb{E}\left[L_g^2||v_T - v^*_T||_2^2 + TL_b^2\int_0^T||v_s - v^*_s||_2^2\,ds\right] \\
             &\leq 8(L^2_g + T^2L_b^2)\mathbb{E}\left[\sup_{0\leq t\leq T}||v_t - v^*_t||_2^2\right] \\
             &\leq 8T(L^2_g + T^2L_b^2)\mathbb{E}\left[\int_0^T||\bar{\rho}(t,\lambda_t) - \bar{\rho}(t,\lambda^*_t)||_2^2\,dt\right]
        \end{align*}
        Now applying \cref{lemma:d2_bound_control} we have:
        \begin{align*}
             \mathbb{E}\left[\sup_{0\leq t\leq T}||\Phi(\lambda)_t - \lambda^*_t||_2^2\right] \leq \frac{8T(L^2_g + T^2L_b^2)}{\mu^2}||\Phi(\lambda) - \lambda||_{\mathbb{L}^2_{\mathbb{F}}}^2
        \end{align*}
        and the result follows.
    \end{proof}    
\end{lemma}
\paragraph{Proof of \cref{lemma:fixed_point_stability}}
    \begin{proof}[Proof of 1]
        By the triangle inequality we have that:
        \begin{align*}
            d_2(\lambda, \lambda^*) \leq d_2(\Phi(\lambda), \lambda) + d_2(\Phi(\lambda), \lambda^*), \quad \lambda \in \mathcal{A}
        \end{align*}
         By \cref{lemma:fixed_point_stability_bkwd_dS} we have that there exists $C > 0$:
    \begin{align*}
        d_2(\Phi(\lambda), \lambda^*) \leq \sqrt{T}d_S(\Phi(\lambda), \lambda^*)\leq \sqrt{T}Cd_2(\Phi(\lambda), \lambda), \quad \forall \lambda \in \mathcal{A}
    \end{align*}
    thus we have that:
    \begin{align*}
         d_2(\lambda, \lambda^*) \leq (1+\sqrt{T}C)d_2(\Phi(\lambda), \lambda), \quad \forall \lambda \in \mathcal{A}
    \end{align*}
    \end{proof}
    \begin{proof}[Proof of 2]
        After fixing c\`adl\`ag representatives by the triangle inequality we have:
        \begin{align*}
            d_S(\lambda, \lambda^*)  \leq d_S(\Phi(\lambda), \lambda) + d_S(\Phi(\lambda), \lambda^*), \quad \lambda \in \mathcal{A}
        \end{align*}
        By \cref{lemma:fixed_point_stability_bkwd_dS} we have that there exists $C > 0$:
        \begin{align*}
            d_S(\Phi(\lambda), \lambda^*) \leq Cd_2(\Phi(\lambda),\lambda) \leq \sqrt{T}Cd_S(\Phi(\lambda), \lambda), \quad \forall \lambda \in \mathcal{A}
        \end{align*}
        thus we have:
        \begin{align*}
            d_S(\lambda, \lambda^*) \leq (1+\sqrt{T}C)d_S(\Phi(\lambda), \lambda), \quad \forall \lambda \in \mathcal{A}
        \end{align*}
    \end{proof}

\subsection{Solvers}
\subsubsection{Martingale Solver}\label{subsec:martingale_solver}

Before proving \cref{theorem:martingale_solver_stable}, we prove that the forward dynamics in the variational problem of \ref{eq:martingale_solver} are well posed for $M \in \mathbb{M}^{2}_{\mathbb{F}}([0,T]; \mathbb{R}^d)$ (in-fact we can prove it on the larger space $\mathbb{S}^{2}_{\mathbb{F}}([0,T]; \mathbb{R}^d)$).
\begin{lemma}[Existence and Uniqueness for Forward System] \label{theorem:approx_forward_eqn_ex_un}
    For every $M \in \mathbb{S}^2_{\mathbb{F}}([0,T]; \mathbb{R}^d)$ there exists a unique solution pair $(v, \lambda) \in \mathbb{S}^{2,c}_{\mathbb{F}}([0,T]; \mathbb{R}^d)\times \mathbb{S}^{2}_{\mathbb{F}}([0,T]; \mathbb{R}^d)$ to the equation governing the forward dynamics in \ref{eq:martingale_solver}.
    \begin{proof}
        To show uniqueness, taking $M \in \mathbb{S}^2_{\mathbb{F}}([0,T]; \mathbb{R}^d)$ let $(v^1, \lambda^1)$ and $(v^2, \lambda^2)$ be solutions and let $\Delta\lambda = \lambda - \lambda$ and $\Delta v = v^1 - v^2$. We have:
        \begin{align*}
            ||\Delta v_t||_2^2 + ||\Delta\lambda_t||_2^2 &= 2\int_0^t\langle \Delta v_s, \bar{\rho}(s,\lambda^1_s) - \bar{\rho}(s,\lambda^2_s) \rangle\,ds - 2\int_0^t\langle \Delta\lambda_s, b(s, v^1_s)- b(s, v^2_s)\rangle\,ds \\
            &\leq 2\int_0^tL_{\rho}||\Delta v_s||_2||\Delta\lambda_s||_2\,ds + 2\int_0^tL_b|| \Delta\lambda_s||_2||\Delta v_s||_2\,ds \\
            &\leq \left(L_{\rho} + L_b\right)\int_0^t\left(||\Delta v_s||_2^2 + ||\Delta\lambda_s||_2^2\right)\,ds, \quad t \in [0,T] 
        \end{align*}
        with:
        \begin{align*}
            \left(L_{\rho} + L_b\right)\int_0^T\left(||\Delta v_s||_2^2 + ||\Delta\lambda_s||_2^2\right)\,ds \leq T\left(L_{\rho} + L_b\right)\left(\sup_{t \in [0,T]}||\Delta v_s||_2^2 + \sup_{0\leq t\leq T}||\lambda_s||_2^2\right)
        \end{align*}
         by assumption, $\Delta v \in \mathbb{S}^{2,c}_{\mathbb{F}}([0,T]; \mathbb{R}^d)$ and $\Delta\lambda \in \mathbb{S}^{2}_{\mathbb{F}}([0,T]; \mathbb{R}^d)$, thus there exists $\bar{\Omega}\subset\Omega$ with $\mathbb{P}(\bar{\Omega}) = 1$ such that:
        \begin{align*}
            T\left(L_{\rho} + L_b\right)\left(\sup_{t \in [0,T]}||\Delta v(s,\omega)||_2^2 + \sup_{0\leq t\leq T}||\lambda(s,\omega)||_2^2\right)< \infty, \quad \forall \omega \in \bar{\Omega}
        \end{align*}
        thus by Gronwall's inequality (\cite{PardouxBook:2014}, Corollary 6.60) we have that:
        \begin{align*}
            ||\Delta v(t,\omega)||_2^2 + ||\Delta\lambda(t,\omega)||_2^2 = 0, \quad \forall t \in [0,T], \quad \forall \omega \in \bar{\Omega}
        \end{align*}
        implying that:
        \begin{align*}
            \mathbb{E}\left[\sup_{t \in [0,T]}(||\Delta v_t||_2^2  + ||\Delta \lambda_t||_2^2)\right]^{1/2} = 0
        \end{align*}
        and we can conclude that:
        \begin{align*}
            \mathbb{P}(v_t^1 = v^2_t, \lambda_t^1 = \lambda_t^2; \forall t \in [0,T]) = 1
        \end{align*}
        proving uniqueness up to indistinguishability. To prove existence, fixing $t \in [0,T]$ and $\delta > 0$ with $t+\delta \leq T$, define the Picard iteration:
        \begin{align*}
            \lambda^{k+1}_s &= \tilde{\lambda} -\int_t^sb(u, v^k_u)\,du +  M_{s} - M_t, \quad s \in [t,t+\delta] \\
            v^{k+1}_s &= \tilde{v} + \int_t^s\bar{\rho}(u, \lambda^{k}_u)\,du, \quad s \in [t,t+\delta]
        \end{align*}
        with $(v^0_s, \lambda^0_s)_{s \in [t,t+\delta]} \equiv(\tilde{v},\tilde{\lambda}) \in L^2(\Omega, \mathcal{F}_t, \mathbb{P}; \mathbb{R}^d) \times L^2(\Omega, \mathcal{F}_t, \mathbb{P}; \mathbb{R}^d)$.  We have:
        \begin{align*}
            ||\Delta \lambda^{k+1}_s||_2 &\leq L_b\int_t^{t+\delta}||\Delta v^k_u||_2\,du \\
            ||\Delta v^{k+1}_s||_2 &\leq L_\rho\int_t^{t+\delta} ||\Delta \lambda^k_u||_2\,du
        \end{align*} 
        thus:
        \begin{align*}
            ||\Delta v^{k+1}_s||_2^2 + ||\Delta \lambda^{k+1}_s||_2^2 \leq \delta^2(L_b + L_\rho)^2\sup_{s \in [t,t+\delta]}\left(||\Delta v^k_u||_2^2 + ||\Delta \lambda^k_u||_2^2\right)
        \end{align*}
        and we have
        \begin{align*}
            \mathbb{E}\left[\sup_{s \in [t,t+\delta]}(||\Delta v^{k+1}_s||_2^2 + ||\Delta \lambda^{k+1}_s||_2^2)\right]^{1/2} \leq \delta(L_b + L_\rho)\mathbb{E}\left[\sup_{s \in [t,t+\delta]}(||\Delta v^{k}_s||_2^2 + ||\Delta \lambda^{k}_s||_2^2)\right]^{1/2}
        \end{align*}
        Choosing $\delta < (L_b + L_\rho)$ we obtain a contraction and the sequence $(v^k, \lambda^k)$ is Cauchy and thus converges to a unique limit $(v,\lambda) \in \mathbb{S}^{2,c}_{\mathbb{F}}([0,T]; \mathbb{R}^d)\times \mathbb{S}^{2}_{\mathbb{F}}([0,T]; \mathbb{R}^d)$. By Banach's fixed point theorem we have:
        \begin{align*}
            \lambda_s &= \tilde{\lambda} -\int_t^sb(u, v_u)\,du +  M_{s} - M_t, \quad s \in [t,t+\delta] \\
            v_s &= \tilde{v} + \int_t^s\bar{\rho}(u, \lambda_u)\,du, \quad s \in [t,t+\delta]
        \end{align*}
        The result follows by starting at $t=0$ with $\tilde{\lambda} = M_0$ and $\tilde{v} = v_0$ we obtaining the solution by applying the procedure iteratively on $[0, \delta_0], [\delta_1, \delta_2], \dots $.
    \end{proof}
\end{lemma}
\paragraph{Proof of \cref{theorem:martingale_solver_stable}}
    \begin{proof}[Proof of 1]
        Let $M \in \mathbb{M}^{2}_{\mathbb{F}}([0,T]; \mathbb{R}^d)$. Solving the forward equation for the pair $(v, \lambda)$ define $M' \in \mathcal{A}$ such that:
        \begin{align*}
            M_t' := \Phi(\lambda)_t - \lambda_t + M_t, \quad t \in [0,T]
        \end{align*}
        Then we have $M'\in \mathbb{M}^{2}_{\mathbb{F}}([0,T];\mathbb{R}^d)$ with:
        \begin{align*}
            M_t' = \mathbb{E}\left[g(v_T) + \int_t^Tb(s,v_s)\,ds \ \bigg|\mathcal{F}_t\right], \quad \mathbb{P}\text{-a.s.}, \quad t \in [0,T]
        \end{align*}
        By Doob's inequality:
        \begin{align*}
            J(M) &= \mathbb{E}[||g(v_T) - \lambda_T||_2^2] \\
            &= \mathbb{E}\left[||M_T' - M_T||_2^2\right] \\
            &\geq \frac{1}{4}\mathbb{E}\left[\sup_{0\leq t\leq T}||M'_t - M_t||\right] \\  
            &= \frac{1}{4}\mathbb{E}\left[\sup_{0\leq t\leq T}||\Phi(\lambda)_t - \lambda_t||\right]  \geq \frac{1}{4T}\mathbb{E}\left[\int_0^T||\Phi(\lambda)_t - \lambda_t||_2^2\right]
        \end{align*}
        Thus we have that:
        \begin{align*}
            d_S(\Phi(\lambda), \lambda) \leq 2J(M)^{1/2}, \quad d_2(\Phi(\lambda), \lambda) \leq 2\sqrt{T}J(M)^{1/2}, \quad \forall M \in \mathbb{M}^{2}_{\mathbb{F}}([0,T]; \mathbb{R}^d)
        \end{align*}
        Choosing $M^* \in \mathbb{M}^{2}_{\mathbb{F}}([0,T]; \mathbb{R}^d)$ such that when solving the forward equation for $(v^*, \lambda^*)$ we have:
        \begin{align*}
            M^*_t &= \mathbb{E}\left[g(v^*_T) + \int_0^Tb(s,v^*_s)\,ds \ \bigg|\mathcal{F}_t\right], \quad \mathbb{P}\text{-a.s.}, \quad t \in [0,T] \\
            \lambda^*_t &= - \int_0^tb(s,v^*_s)\,ds + M^*_t \\
            v^*_t &= v_0 + \int_0^t\bar{\rho}(s,\lambda^*_s)\,ds
        \end{align*}
        Then clearly we have that $J(M^*) = 0$, and $(v^*,\lambda^*)$ solves \cref{eq:fbs_general}. The result follows from non-negativity of $J$ and existence and uniqueness of such an $M^*$ which follows from Theorem 5.5 and Theorem 5.11 in \cite{KeyworthPaper:2026}.    
    \end{proof}
    \begin{proof}[Proof of 2]
        We have by 1. that      
        \begin{align*}
            d_S(\Phi(\lambda), \lambda) \leq 2J(M)^{1/2}, \quad d_2(\Phi(\lambda), \lambda) \leq 2\sqrt{T}J(M)^{1/2} \quad \forall M \in \mathbb{M}^{2}_{\mathbb{F}}([0,T]; \mathbb{R}^d)
        \end{align*}
    By \cref{lemma:fixed_point_stability} we have that there exists $C > 0$ such that:
    \begin{align*}
            d_S(\lambda, \lambda^*) 
            \leq Cd_S(\Phi(\lambda), \lambda), \quad d_2(\lambda, \lambda^*) 
            \leq Cd_2(\Phi(\lambda), \lambda), \quad \forall \lambda \in \mathcal{A} 
    \end{align*}
    Thus we have that:
     \begin{align*}
        d_S(\lambda, \lambda^*) \leq 2CJ(M)^{1/2}, \quad d_2(\lambda, \lambda^*) \leq 2\sqrt{T}CJ(M)^{1/2} \quad \forall M \in \mathbb{M}^{2}_{\mathbb{F}}([0,T]; \mathbb{R}^d)
    \end{align*}
    \end{proof}    
\subsubsection{Hilbert Space Orthogonal Min-Max Solvers} \label{subsec:min_max_solver}
Before proving \cref{theorem:unique_sol_const_opt} and \cref{theorem:min_max_game} we prove a useful lemma.
\begin{lemma} \label{lemma:proj_target_mart}
    For $\lambda \in \mathcal{A}$, $\pi Y \in \mathcal{A}$ is a square integrable martingale such that:
    \begin{align*}
            \pi Y_t = \mathbb{E}\left[g(v_T) + \int_0^Tb(s, v_s)\,ds\bigg|\mathcal{F}_t\right], \quad \mathbb{P}\text{-a.s.}, \  t \in [0,T]
    \end{align*}
    \begin{proof}
        By \cref{lemma:orth_proj} we have that:
        \begin{align*}
            \pi Y_t = \mathbb{E}[Y_t|\mathcal{F}_t]= \mathbb{E}\left[g(v_T) + \int_0^Tb(s, v_s)\,ds\bigg|\mathcal{F}_t\right], \quad \mathbb{P}\text{-a.s.}, \ \text{a.e. } t \in [0,T]
        \end{align*}
        The RHS is a square integrable martingale that admits a right continuous version (which exists by the corollary of Theorem 3.1 in \cite{LipsterBook1:2001}), thus by the definition of $\mathcal{A}$ we have that $\pi Y \in \mathcal{A}$. 
    \end{proof}
\end{lemma}

\paragraph{Proof of \cref{theorem:unique_sol_const_opt}}
\begin{proof}
    Suppose for $i \in \{1,2\}$ $\lambda^i \in \mathcal{A}$ satisfy the constraint $F(\lambda^i) = 0$, then by \cref{lemma:orth_proj}:
    \begin{align*}
        \pi Y^i_t = M_t^i, \quad \mathbb{P}\text{-a.s.}, \ \text{a.e. } t \in [0,T], \quad i \in \{1,2\}
    \end{align*}
     Since $M^i\in \mathcal{A}$ and since by \cref{lemma:proj_target_mart} $\pi Y^i \in \mathcal{A}$, we can conclude by standard arguments that $M^i$ is indistinguishable from a right continuous version of the square integrable martingale:
    \begin{align*}
        M^{i}_t = \mathbb{E}\left[g(v^{i}_T) + \int_0^Tb(s,v^{i}_s)\,ds \ \bigg|\mathcal{F}_t\right], \quad \mathbb{P}\text{-a.s.}, \quad  t \in [0,T]
    \end{align*}
    since we have forced the dynamics, by \ref{eq:controlled_dynamics} we have that $(\lambda^i, v^i)$ solve \cref{eq:fbs_general}: 
    \begin{align*}
        v_t^{i} &= v_0 + \int_0^t\bar{\rho}(s, \lambda^i_s)\,ds, \quad t \in [0,T] \\
        \lambda_t^i &= \mathbb{E}\left[g(v^{i}_T) + \int_t^Tb(s, v^{i}_s)\,ds \ \bigg|\mathcal{F}_t\right], \quad \mathbb{P}\text{-a.s.}, \ t \in [0,T]
    \end{align*}
    By non-negativity of $J$ and the terminal condition for $\lambda^i$ we have that:
    \begin{align*}
        \inf_{\lambda \in \mathcal{A}}J(\lambda) = J(\lambda^i) = 0, \quad i \in \{1, 2\}
    \end{align*}
    and thus $\lambda^1, \lambda^2$ are solutions of the constrained optimisation of \ref{eq:const_opt}, but by Theorem 5.5 in \cite{KeyworthPaper:2026} we must have that:
    \begin{align*}
        \mathbb{P}(v^{1}_t = v^{2}_t; \forall t \in [0,T]) = 1, \quad \mathbb{P}(\lambda^1_t = \lambda^2_t; \forall t \in [0,T]) = 1
    \end{align*}
    proving uniqueness. Existence follows from Theorem 5.11 in \cite{KeyworthPaper:2026}.
\end{proof}
\paragraph{Proof of \cref{theorem:min_max_game}}
    \begin{proof}[Proof of 1]
        By \cref{theorem:unique_sol_const_opt} there exists a unique $\lambda^* \in \mathcal{A}$ such that:
        \begin{align*}
            J(\lambda^*) = 0, \quad F(\lambda^*) = 0
        \end{align*}
        Furthermore $\lambda^* \in \operatorname{Fix}\Phi$ and the pair $(v^*, \lambda^*)$ solves \cref{eq:fbs_general}. By \cref{corr:orth_proj_cost} we can conclude that:
        \begin{align*}
            F(\lambda^*) = 0 \iff \sup_{\varphi \in \mathcal{T}_c}\langle \varphi, Y^* - M^{*}\rangle_{\mathbb{L}^2} = 0
        \end{align*}
         Thus we have $\lambda^*$ uniquely satisfies:
        \begin{align*}
            \sup_{\varphi \in \mathcal{T}_c}L(\lambda^*, \varphi) = J(\lambda^*) + \sup_{\varphi \in \mathcal{T}_c}\langle \varphi, Y^* - M^*\rangle_{\mathbb{L}^2} = 0
        \end{align*}
         Suppose now that for some $\lambda \in \mathcal{A}$ we have:
        \begin{align*}
            \sup_{\varphi \in \mathcal{T}_c}L(\lambda, \varphi) < \sup_{\varphi \in \mathcal{T}_c}L(\lambda^*, \varphi) = 0
        \end{align*}
        Then we must have that (since $J$ is non-negative):
        \begin{align*}
            \sup_{\varphi \in \mathcal{T}_c}\langle \varphi, Y - M\rangle_{\mathbb{L}^2} < -J(\lambda) \leq 0
        \end{align*}
        But we have that:
        \begin{align*}
            \sup_{\varphi \in \mathcal{T}_c}\langle \varphi, Y - M\rangle_{\mathbb{L}^2} = \sup_{\varphi \in \mathcal{T}_c}\langle \varphi, \pi Y - M\rangle_{\mathbb{L}^2_{\mathbb{F}}} = c||\pi Y - M||_{\mathbb{L}^2_{\mathbb{F}}} = c||\Phi(\lambda) - \lambda||_{\mathbb{L}^2_{\mathbb{F}}} \geq 0
        \end{align*}
        a contradiction. Thus we have that:
        \begin{align*}
            0 = \sup_{\varphi \in \mathcal{T}_c}L(\lambda^*, \varphi) = \inf_{\lambda \in \mathcal{A}}\sup_{\varphi \in \mathcal{T}_c}L(\lambda, \varphi) 
        \end{align*}
        and, with $C_1 = 1/c$, we have that:
        \begin{align*}
            d_2(\Phi(\lambda),\lambda) \leq C_1\sup_{\varphi \in \mathcal{T}_c}L(\lambda, \varphi), \quad \forall \lambda \in \mathcal{A}
        \end{align*}
        Now, for $\lambda \in \mathcal{A}$ we have:
        \begin{align*}
            \sup_{\varphi \in \mathcal{D}}\tilde{L}(\lambda, \varphi) &= J(\lambda) + \sup_{\varphi \in \mathcal{D}}\left(\langle \varphi, Y - M\rangle_{\mathbb{L}^2} - \frac{\tilde{c}}{2}||\varphi||_{\mathbb{L}^2_{\mathbb{F}}}^2\right) \\
            &= J(\lambda) + \sup_{\varphi \in \mathcal{D}}\left(\langle \varphi, \pi Y - M\rangle_{\mathbb{L}^2_{\mathbb{F}}} - \frac{\tilde{c}}{2}||\varphi||_{\mathbb{L}^2_{\mathbb{F}}}^2\right) \\
            &= J(\lambda) + \langle \varphi^*, \pi Y - M\rangle_{\mathbb{L}^2_{\mathbb{F}}} - \frac{\tilde{c}}{2}||\varphi^*||_{\mathbb{L}^2_{\mathbb{F}}}^2
        \end{align*}
        Where we have by standard arguments that $\varphi^* = (\pi Y - M)/r$. Thus by \cref{lemma:proj_target_mart} and \ref{eq:controlled_dynamics}:
        \begin{align*}
            \sup_{\varphi \in \mathcal{D}}\tilde{L}(\lambda, \varphi) = J(\lambda) + \frac{1}{2r}||\pi Y - M||_{\mathbb{L}^2_{\mathbb{F}}}^2 = J(\lambda) + \frac{1}{2r}||\Phi(\lambda) - \lambda||_{\mathbb{L}^2_{\mathbb{F}}}^2, \quad \forall \lambda \in \mathcal{A}
        \end{align*}
        Now taking the infimum, by non-negativity and existence and uniqueness of $\lambda^*$ we can conclude that:
        \begin{align*}
            0 = J(\lambda^*) + \frac{1}{2r}||\Phi(\lambda^*) - \lambda^*||_{\mathbb{L}^2_{\mathbb{F}}}^2 = \sup_{\varphi \in \mathcal{D}}\tilde{L}(\lambda^*, \varphi) = \inf_{\lambda \in \mathcal{A}}\sup_{\varphi \in \mathcal{D}}\tilde{L}(\lambda, \varphi)
        \end{align*}
        we also have, with $C_2 = 2r$, that:
        \begin{align*}
            d_2(\Phi(\lambda),\lambda)^2 \leq C_2\sup_{\varphi \in \mathcal{D}}\tilde{L}(\lambda, \varphi), \quad \forall \lambda \in \mathcal{A}
        \end{align*}
    \end{proof}
    \begin{proof}[Proof of 2]
        By 1. we have that there exist $C_1, C_2 > 0$ such that:
        \begin{align*}
            d_2(\Phi(\lambda), \lambda) \leq C_1\sup_{\varphi \in \mathcal{T}_c}L(\lambda, \varphi), \quad d_2(\Phi(\lambda), \lambda) \leq C_2\sup_{\varphi \in \mathcal{D}}\tilde{L}(\lambda, \varphi), \quad \forall \lambda \in \mathcal{A}
        \end{align*}
        By \cref{lemma:fixed_point_stability} there exists $C > 0$ such that:
        \begin{align*}
            d_2(\lambda, \lambda^*) 
            \leq Cd_2(\Phi(\lambda), \lambda), \quad \forall \lambda \in \mathcal{A}
        \end{align*}
        and the result follows.
    \end{proof}
\subsection{Discretized Solvers} \label{subsec:discretized_solvers}
First some useful lemmas that will help us make the connection between the objective functions of the discretized solvers and the fixed point stability results.
\begin{remark}\label{remark:tilde_m}
      For $\lambda^N \in \mathcal{A}$ we have $\lambda^N \in \mathcal{A}$ recall that $\Phi:\mathcal{A}\rightarrow \mathcal{A}$ so $\Phi(\lambda^N) \in \mathcal{A}$ and for $t \in [0,T]$ we have:
    \begin{align*}
        \Phi(\lambda^N)_t &= \mathbb{E}\left[g(v_T) + \int_t^Tb(s,v_s)\,ds \ \bigg|\mathcal{F}_t\right], \quad \mathbb{P}\text{-a.s.}\\
        v_t &= v_0 + \int_0^t\mathcal{N}(\lambda^N)_s\,ds = v_0 + \int_0^t\bar{\rho}(s,\lambda^N_s)\,ds
    \end{align*}
    we can also define martingales $\tilde{M}, \tilde{M}^N \in \mathbb{M}^{2}_{\mathbb{F}}([0,T]; \mathbb{R}^d)$ such that:
    \begin{align*}
        \tilde{M}_t &= \mathbb{E}\left[g(v_T) + \int_0^Tb(s,v_s)\,ds \ \bigg|\mathcal{F}_t\right], \quad \mathbb{P}\text{-a.s.}, \quad t \in [0,T] \\
        \tilde{M}_t^N &= \mathbb{E}\left[g(v_{t_N}^N) + \int_0^Tb(\tau(s), v_{\tau(s)}^N)\,ds \ \bigg|\mathcal{F}_t\right], \quad \mathbb{P}\text{-a.s.}, \quad t \in [0,T]
    \end{align*}
\end{remark}
\begin{lemma}\label{lemma:discrete_fpr}
     If for $M^N, \lambda^N \in \mathcal{A}$ the following:
    \begin{align*}
        \lambda_t^N &= -\int_0^tb(\tau(s), v_{\tau(s)}^N)\,ds + M_t^N, \quad t \in [0,T] \\
        v_t^N &= v_0 + \int_0^t\bar{\rho}(\tau(s), \lambda_{\tau(s)}^N)\,ds, \quad t \in [0,T]
    \end{align*}    
    then $\Phi(\lambda^N) - \lambda^N \in \mathcal{A}$ with:
    \begin{align*}
        \Phi(\lambda^N)_t - \lambda^N_t = \mathbb{E}\left[g(v_T) - g(v_{t_N}^N) + \int_t^T(b(s, v_s) - b(\tau(s),v_{\tau(s)}^N)\,ds \ \bigg| \mathcal{F}_t\right] \\ {+ (\tilde{M}_t^N - M_t^N), \quad \mathbb{P}\text{-a.s.}, \ t \in [0,T]}
    \end{align*}
    and:
    \begin{align*}
        \mathbb{E}\left[\sup_{0\leq t\leq T}||\mathbb{E}\left[g(v_T) - g(v_{t_N}^N) + \int_t^T(b(s, v_s) - b(\tau(s),v_{\tau(s)}^N))\,ds \ \bigg| \mathcal{F}_t\right]||_2^2 \right]\\ 
        \leq 8\mathbb{E}\left[||g(v_T) - g(v_{t_N}^N)||^2_2 + T\int_0^T||b(s, v_s) - b(\tau(s),v_{\tau(s)}^N)||_2^2\,ds\right]    
    \end{align*}
    \begin{proof}
        Clearly with $\tilde{M}$ as in \cref{remark:tilde_m} we have that:
        \begin{align*}
            \Phi(\lambda^N)_t = - \int_0^tb(s,v_s)\,ds + \tilde{M}_t, \quad \mathbb{P}\text{-a.s.}, \quad t \in [0,T]
        \end{align*}
        thus we have $\mathbb{P}\text{-a.s.}, \ t \in [0,T]$  that:
        \begin{align*}
            \Phi(\lambda^N)_t - \lambda^N_t &= - \int_0^t( b(s,v_s)-b(\tau(s), v_{\tau(s)}^N))\,ds + \tilde{M}_t - M_t^N \\
            &= - \int_0^t( b(s,v_s)-b(\tau(s), v_{\tau(s)}^N))\,ds + (\tilde{M}_t - \tilde{M}^N_t) + (\tilde{M}^N_t - M_t^N)
        \end{align*}
        hence:
        \begin{align*}
            \Phi(\lambda^N)_t - \lambda^N_t = \mathbb{E}\left[g(v_T) - g(v_{t_N}^N) + \int_t^T(b(s, v_s) - b(\tau(s),v_{\tau(s)}^N)\,ds \ \bigg| \mathcal{F}_t\right] \\ {+ (\tilde{M}_t^N - M_t^N), \quad \mathbb{P}\text{-a.s.}, \ t \in [0,T]}
        \end{align*}
        Consequently we obtain:
        \begin{align*}
            &\sup_{0\leq t\leq T}||\Phi(\lambda^N)_t - \lambda^N_t - (\tilde{M}_t^N - M_t^N)||_2^2 \\ &\leq \sup_{0\leq t\leq T}\mathbb{E}\left[||g(v_T) - g(v_{t_N}^N)|| + \int_0^T||b(s, v_s) - b(\tau(s),v_{\tau(s)}^N)||_2\,ds \ \bigg| \mathcal{F}_t\right]^2, \quad \mathbb{P}\text{.a.s}, \quad t \in [0,T]
        \end{align*}
        and by the same argument as in the proof of \cref{lemma:bounding_mart} by Doob's submartingale inequality we obtain
        \begin{align*}
            &\mathbb{E}\left[\sup_{0\leq t\leq T}||\Phi(\lambda^N)_t - \lambda^N_t - (\tilde{M}_t^N - M_t^N)||_2^2\right] \\
            &\quad\quad \leq 8\mathbb{E}\left[||g(v_T) - g(v_{t_N}^N)||^2_2 + T\int_0^T||b(s, v_s) - b(\tau(s),v_{\tau(s)}^N)||_2^2\,ds\right]
        \end{align*}
    \end{proof}
\end{lemma}
\subsubsection{Martingale Solver}\label{subsec:martingale_solver_convergence}
We now prove \cref{theorem:martingale_solver_discrete}.
\paragraph{Proof of \cref{theorem:martingale_solver_discrete}}
    \begin{proof}
        We have that $\forall M^N \in \mathbb{M}^{2}_{\mathbb{F}}([0,T]; \mathbb{R}^d)$:
        \begin{align*}
            g(v_{t_N}^N) - \lambda_{t_N}^N = g(v_{t_N}^N) + \int_0^Tb(\tau(s), v_{\tau(s)}^N)\,ds - M_{t_N}^N = \tilde{M}_{t_N}^N - M_{t_N}^N
        \end{align*}
        With $\tilde{M}^N$ defined as in \cref{remark:tilde_m} and thus $\tilde{M}^N - M^N \in \mathbb{M}^{2}_{\mathbb{F}}([0,T]; \mathbb{R}^d)$. By Doob's inequality:
        \begin{align*}
            J_N(M^N) = \mathbb{E}\left[||g(v_{t_N}^N) - \lambda_{t_N}^N||_2^2\right] = \mathbb{E}\left[||\tilde{M}_{t_N}^N - M_{t_N}^N||_2^2\right] \geq \frac{1}{4}\mathbb{E}\left[\sup_{0\leq t\leq T}||\tilde{M}_{t}^N - M_{t}^N||_2^2\right] 
        \end{align*}
        Now defining:
        \begin{align*}
            \delta_N := \sup_{M^N \in \mathbb{M}^{2}_{\mathbb{F}}}\mathbb{E}\left[||g(v_T) - g(v_{t_N}^N)||^2_2 + T\int_0^T||b(s, v_s) - b(\tau(s),v_{\tau(s)}^N)||_2^2\,ds\right]
        \end{align*}
        It follows from \cref{lemma:discrete_fpr} that:
        \begin{align*}
            &\mathbb{E}\left[\sup_{0\leq t\leq T}||\Phi(\lambda^N)_t - \lambda^N_t ||_2^2\right] \\
             &\quad\quad \leq 2\mathbb{E}\left[\sup_{0\leq t\leq T}\bigg|\bigg|\mathbb{E}\left[g(v_T) - g(v_{t_N}^N) + \int_t^T(b(s, v_s) - b(\tau(s),v_{\tau(s)}^N)\,ds \ \bigg| \mathcal{F}_t\right]\bigg|\bigg|_2^2\right] \\ 
             &\quad\quad \leq 2\mathbb{E}\left[\sup_{0\leq t\leq T}||\tilde{M}_t^N - M_t^N||_2^2\right] \leq 8(J_N(M^N) + 2\delta_N)
        \end{align*}
        By \cref{lemma:fixed_point_stability} we have that there exists $C > 0$ such that $d_S(\lambda, \lambda^*)^2 
        \leq Cd_S(\Phi(\lambda), \lambda)^2, \ \forall \lambda \in \mathcal{A}$ and we obtain:
        \begin{align*}
            d_S(\lambda^N, \lambda^*)^2 \leq 8C(J_N(M^N) + 2\delta_N), \quad \forall M^N \in \mathbb{M}^{2}_{\mathbb{F}}([0,T]; \mathbb{R}^d)
        \end{align*}
        Now we have that:
        \begin{align*}
            J_N(M^N) = \mathbb{E}\left[||\tilde{M}_{t_N}^N - M_{t_N}^N||_2^2\right] \leq d_S(\tilde{M}^N, M^N)^2, \quad \forall M^N \in \mathscr{A}_M\subseteq \mathbb{M}^{2}_{\mathbb{F}}([0,T]; \mathbb{R}^d)
        \end{align*}
        Thus we have:
        \begin{align*}
                0\leq \inf_{M^N \in \mathscr{A}_M}J_N(M^N) \leq  \inf_{M^N \in \mathscr{A}_M}d_S(\tilde{M}^N, M^N)^2 \leq \sup_{M \in \mathbb{M}^{2}_{\mathbb{F}}}d_S(M, \mathscr{A}_M)^2
        \end{align*}
         where we have:
        \begin{align*}
            d_S(M, \mathscr{A}_M) := \inf_{\tilde{M} \in \mathscr{A}_M}d_S(M, \tilde{M}), \quad  M \in \mathbb{M}^{2}_{\mathbb{F}}([0,T]; \mathbb{R}^d)
        \end{align*}
    \end{proof}

\subsubsection{Discretized Min-Max Solver}\label{subsec:min_max_solver_convergence}
Before proving \cref{theorem:discretized_min_max_solvers} we prove an important lemma.
\begin{lemma}\label{lemma:bounded_Y_M_min_max_solver}
    Let $C_\lambda > 0$.
    \begin{enumerate}
        \item Then with $Y^N, M^N$ as in \ref{eq:discretized_dynamics_min_max_solver} there exist constants $C_1, C_2 > 0$ such that:
        \begin{align*}
            ||Y^N||_{\mathbb{L}^2_{\mathbb{F}}} \leq C_1, \quad ||M^N||_{\mathbb{L}^2_{\mathbb{F}}} \leq C_2, \quad \forall ||\lambda^N||_{\mathbb{S}^2_{\mathbb{F}}} \leq C_{\lambda}
        \end{align*}
        \item Then there exists $C > 0$ such that:
        \begin{align*}
            ||\Phi(\lambda^N)||_{\mathbb{L}^2_{\mathbb{F}}} \leq C, \quad \forall ||\lambda^N||_{\mathbb{S}^2_{\mathbb{F}}} \leq C_{\lambda}
        \end{align*}
    \end{enumerate}
    \begin{proof}[Proof of 1]
        Setting $v^N$ as in \ref{eq:discretized_dynamics_min_max_solver} with $v_0 = x \in \mathbb{R}^d$, then $s \in [0,T]$ we have:
        \begin{align*}
            || v_{s}^N||_2^2 &\leq 2\left(||x||_2^2 + T\int_0^T||\bar{\rho}(\tau(s),\lambda^N_{\tau(s)})||_2^2\,ds\right) \\
            &\leq 2||x||_2^2 + 4T\int_0^T|| \bar{\rho}(\tau(s),0)||_2^2\,ds + 4TL_\rho^2\int_0^T||\lambda^N_{\tau(s)}||_2^2\,ds
        \end{align*}
        and by \ref{ass:I1} that:
        \begin{align*}
            \mathbb{E}\left[\sup_{s \in [0,T]}|| v_{s}^N||_2^2\right] \leq 2||x||_2^2 + 4T^2\mathbb{E}\left[\sup_{s \in [0,T]}||\bar{\rho}(s,0)||_2^2\right] + 4T^2L_\rho^2C_{\lambda} < \infty, \quad \forall x \in \mathbb{R}^d
        \end{align*}
        Now setting:
        \begin{align*}
            C_{x}= \left(2||x||_2^2 + 4T^2\mathbb{E}\left[\sup_{s \in [0,T]}||\bar{\rho}(s,0)||_2^2\right] + 4T^2L_\rho^2C_{\lambda} \right)^{1/2}, \quad x \in \mathbb{R}^d
        \end{align*}
        we have:
        \begin{align*}
            ||Y^N||_{\mathbb{L}^2_{\mathbb{F}}} \leq \sqrt{T}\mathbb{E}\left[||g(v_{t_N}^N)||_2^2\right]^{1/2} + T\mathbb{E}\left[\int_0^T||b(\tau(s), v_{\tau(s)}^N)||_2^2\,ds\right]^{1/2} 
        \end{align*}
        By \ref{ass:A1} and \ref{ass:A5}, we have:
        \begin{align*}
            \mathbb{E}\left[||g(v_{t_N}^N)||_2^2\right]^{1/2} &\leq \mathbb{E}\left[||g(v_0)||_2^2\right]^{1/2} + L_g\mathbb{E}\left[\sup_{s \in [0,T]}||v_{s}^N - v_0||_2^2\right]^{1/2} \\
            &\leq \mathbb{E}\left[||g(v_0)||_2^2\right]^{1/2} + L_gC_{0}  < \infty
        \end{align*}
        by \ref{ass:A5} and \ref{ass:I1} we have:
        \begin{align*}
             \mathbb{E}\left[\int_0^T||b(\tau(s), v_{\tau(s)}^N)||_2^2\,ds\right]^{1/2} &\leq \sqrt{T}\mathbb{E}\left[\sup_{s \in [0,T]}||b(s, 0)||_2^2\right]^{1/2} + \sqrt{T}L_b\mathbb{E}\left[\sup_{s \in [0,T]}|| v_{s}^N||_2^2\right]^{1/2} \\
             &\leq \sqrt{T}\mathbb{E}\left[\sup_{s \in [0,T]}||b(s, 0)||_2^2\right]^{1/2} + \sqrt{T}L_bC_{v_0} < \infty
        \end{align*}
        Consequently we have that:
        \begin{align*}
            ||Y^N||_{\mathbb{L}^2_{\mathbb{F}}} \leq \sqrt{T}\mathbb{E}\left[||g(v_0)||_2^2\right]^{1/2} +  T^{3/2}\mathbb{E}\left[\sup_{s \in [0,T]}||b(s, 0)||_2^2\right]^{1/2} + \sqrt{T}L_gC_{0} + T^{3/2}L_bC_{v_0} < \infty
        \end{align*}
        Now we have by \ref{ass:A5} that:
        \begin{align*}
            ||M^N||_{\mathbb{L}^2_{\mathbb{F}}} \leq ||\lambda^N||_{\mathbb{L}^2_{\mathbb{F}}} + T\mathbb{E}\left[\int_0^T||b(\tau(s),0)||_2^2\,ds\right]^{1/2} + T^{3/2}L_b\mathbb{E}\left[\sup_{s \in [0,T]}||v_{s}^N||_2^2\right]^{1/2}
        \end{align*}
        and by \ref{ass:I1} and previous calculations we have:
        \begin{align*}
            ||M^N||_{\mathbb{L}^2_{\mathbb{F}}} \leq \sqrt{T}C_{\lambda} + T\mathbb{E}\left[\int_0^T||b(\tau(s),0)||_2^2\,ds\right]^{1/2} + T^{3/2}L_bC_{v_0} < \infty
        \end{align*}
    \end{proof}
    \begin{proof}[Proof of 2]
        $\Phi$ is Lipschitz as a composition of Lipschitz maps $\Phi:R^{-1}\circ A\circ \mathcal{N}$ (see \cite{KeyworthPaper:2026}) and the result follows.
    \end{proof}
\end{lemma}
\paragraph{Proof of \cref{theorem:discretized_min_max_solvers}}
    \begin{proof}[Proof of 1]
        Taking $\lambda^N \in \mathcal{A}, \ ||\lambda||_{\mathbb{S}^2_{\mathbb{F}}} \leq C_\lambda$, and defining $(v^N, M^N)$ as in \ref{eq:discretized_dynamics_min_max_solver}, then by \cref{lemma:discrete_fpr} we have that $\Phi(\lambda^N) - \lambda^N \in \mathcal{A}$ with:
         \begin{align*}
            \Phi(\lambda^N)_t - \lambda^N_t = \mathbb{E}\left[g(v_T) - g(v_{t_N}^N) + \int_t^T(b(s, v_s) - b(\tau(s),v_{\tau(s)}^N)\,ds \ \bigg| \mathcal{F}_t\right] \\ {+ (\tilde{M}_t^N - M_t^N), \quad \mathbb{P}\text{-a.s.}, \ t \in [0,T]}
        \end{align*}
        Where $\tilde{M}^N - M^N \in \mathcal{A}$ but not necessarily in $\mathbb{M}^{2}_{\mathbb{F}}([0,T]; \mathbb{R}^d)$. Letting
        \begin{align*}
            X_t = \mathbb{E}\left[g(v_T) - g(v_{t_N}^N) + \int_t^T(b(s, v_s) - b(\tau(s),v_{\tau(s)}^N)\,ds \ \bigg| \mathcal{F}_t\right], \quad t \in [0,T]
        \end{align*}
        and taking $\varphi \in \mathcal{T}_c$ we have that:
        \begin{align*}
            &\mathbb{E}\left[\int_0^T\langle \varphi_t, \Phi(\lambda^N)_t - \lambda^N_t\rangle\,dt\right] = \mathbb{E}\left[\int_0^T\langle \varphi_t, \tilde{M}_t^N - M^N_t\rangle\,dt\right] \\
            &\quad\quad+\mathbb{E}\left[\int_0^T\langle \varphi_t, X_t \rangle\,dt\right]
        \end{align*}
        By Cauchy-Schwarz and \cref{lemma:discrete_fpr} we have:
        \begin{align*}
            \mathbb{E}\left[\int_0^T\langle \varphi_t, X_t \rangle\,dt\right] &\leq ||\varphi||_{\mathbb{L}^2_{\mathbb{F}}}\mathbb{E}\left[\int_0^T ||X_t||_2^2 \,dt\right]^{1/2} \\
            &\leq c\sqrt{T}\mathbb{E}\left[\sup_{0\leq t\leq T}||X_t||_2^2\right]^{1/2} \\
            &\leq  c\sqrt{8T}\mathbb{E}\left[||g(v_T) - g(v_{t_N}^N)||^2_2 + T\int_0^T||b(s, v_s) - b(\tau(s),v_{\tau(s)}^N)||_2^2\,ds\right]^{1/2}
        \end{align*}
        Defining $\delta_N \geq 0$ such that:
        \begin{align*}
            \delta_N := \sup_{\lambda^N \in \mathcal{A}}\mathbb{E}\left[||g(v_T) - g(v_{t_N}^N)||^2_2 + T\int_0^T||b(s, v_s) - b(\tau(s),v_{\tau(s)}^N)||_2^2\,ds\right]
        \end{align*}
        Thus we obtain:
        \begin{align*}
            \mathbb{E}\left[\int_0^T\langle \varphi_t, \Phi(\lambda^N)_t - \lambda^N_t\rangle\,dt\right] \leq  \mathbb{E}\left[\int_0^T\langle \varphi_t, \tilde{M}_t^N - M^N_t\rangle\,dt\right] + c(8T\delta_N)^{1/2}
        \end{align*}
        Now recalling from \cref{remark:tilde_m} that:
        \begin{align*}
            \tilde{M}_t^N &= \mathbb{E}\left[g(v_{t_N}^N) + \int_0^Tb(\tau(s), v_{\tau(s)}^N)\,ds \ \bigg|\mathcal{F}_t\right] = \mathbb{E}[Y^N_t|\mathcal{F}_t], \quad \mathbb{P}\text{-a.s.}, \quad t \in [0,T]
        \end{align*}
        by Tonelli-Fubini and the Tower property we obtain:
        \begin{align*}
            \mathbb{E}\left[\int_0^T\langle \varphi_t, \Phi(\lambda^N)_t - \lambda^N_t\rangle\,dt\right] \leq  \mathbb{E}\left[\int_0^T\langle \varphi_t, Y_t^N - M^N_t\rangle\,dt\right] + c(8T\delta_N)^{1/2}
        \end{align*}
         Now  we have that for every $(\varphi, \tilde{\varphi}) \in \mathcal{T}_c\times \tilde{\mathcal{T}}_c$:
        \begin{align*}
            \mathbb{E}\left[\int_0^T\langle \varphi_t, \Phi(\lambda^N)_t - \lambda^N_t\rangle\,dt\right] \leq  \mathbb{E}\left[\int_0^T\langle \tilde{\varphi}_t, Y_t^N - M^N_t\rangle\,dt\right] \\ + \mathbb{E}\left[\int_0^T\langle \varphi_t-\tilde{\varphi}_t, Y_t^N - M^N_t\rangle\,dt\right] + c(8T\delta_N)^{1/2}
        \end{align*}
        Now since $||\lambda^N||_{\mathbb{S}^2_{\mathbb{F}}} \leq C_{\lambda}$, by Cauchy-Schwarz, the triangle inequality and \cref{lemma:bounded_Y_M_min_max_solver} we have that there exists $C' > 0$ such that:
        \begin{align*}
            \mathbb{E}\left[\int_0^T\langle \varphi_t- \tilde{\varphi}_t, Y_t^N - M^N_t\rangle\,dt\right] &\leq \mathbb{E}\left[\int_0^T||\varphi_t-\tilde{\varphi}_t ||_2^2\,dt\right]^{1/2}\mathbb{E}\left[\int_0^T||Y_t^N - M^N_t||_2^2\,dt\right]^{1/2} \\
            &\leq C'd_2(\varphi, \tilde{\varphi})
        \end{align*}
        Thus $\forall(\varphi, \tilde{\varphi}) \in \mathcal{T}_c\times \tilde{\mathcal{T}}_c$, we have that:
         \begin{align*}
            \mathbb{E}\left[\int_0^T\langle \varphi_t, \Phi(\lambda^N)_t - \lambda^N_t\rangle\,dt\right] \leq  \mathbb{E}\left[\int_0^T\langle \tilde{\varphi}_t, Y_t^N - M^N_t\rangle\,dt\right] + C'd_2(\varphi, \tilde{\varphi}) + c(8T\delta_N)^{1/2}
        \end{align*}
        and thus $\forall(\varphi, \tilde{\varphi}) \in \mathcal{T}_c\times \tilde{\mathcal{T}}_c$:
        \begin{align*}
            \mathbb{E}\left[\int_0^T\langle \varphi_t, \Phi(\lambda^N)_t - \lambda^N_t\rangle\,dt\right] \leq  \sup_{\tilde{\varphi} \in \tilde{\mathcal{T}}_c}\mathbb{E}\left[\int_0^T\langle \tilde{\varphi}_t, Y_t^N - M^N_t\rangle\,dt\right] + C'd_2(\varphi, \tilde{\varphi}) + c(8T\delta_N)^{1/2}
        \end{align*}
        defining $d_2(\varphi, \tilde{\mathcal{T}}_c) := \inf\{d_2(\varphi, \tilde{\varphi}): \tilde{\varphi} \in \tilde{\mathcal{T}}_c\}$ we have that:
        \begin{align*}
            \mathbb{E}\left[\int_0^T\langle \varphi_t, \Phi(\lambda^N)_t - \lambda^N_t\rangle\,dt\right] \leq  \sup_{\tilde{\varphi} \in \tilde{\mathcal{T}}_c}\mathbb{E}\left[\int_0^T\langle \tilde{\varphi}_t, Y_t^N - M^N_t\rangle\,dt\right] + C'd_2(\varphi, \tilde{\mathcal{T}}_c) + c(8T\delta_N)^{1/2}
        \end{align*}
        Thus it follows by standard arguments that:
        \begin{align*}
            cd_2(\Phi(\lambda^N), \lambda^N) &= 
            \sup_{\varphi \in \mathcal{T}_c}\mathbb{E}\left[\int_0^T\langle \varphi_t, \Phi(\lambda^N)_t - \lambda^N_t\rangle\,dt\right] \\
            &\leq  \sup_{\varphi \in \tilde{\mathcal{T}}_c}\mathbb{E}\left[\int_0^T\langle \tilde{\varphi}_t, Y_t^N - M^N_t\rangle\,dt\right] + C'\sup_{\varphi \in \mathcal{T}_c}d_2(\varphi, \tilde{\mathcal{T}}_c) + c(8T\delta_N)^{1/2}
        \end{align*}
        By \cref{lemma:fixed_point_stability} there exists $C > 0$ such that $d_2(\lambda, \lambda^*) \leq Cd_2(\Phi(\lambda), \lambda), \ \forall \lambda \in \mathcal{A}$. Consequently we obtain:
        \begin{align*}
            d_2(\lambda^N, \lambda^*) 
            &\leq  C\left(\frac{1}{c}\sup_{\varphi \in \tilde{\mathcal{T}}_c}\mathbb{E}\left[\int_0^T\langle \tilde{\varphi}_t, Y_t^N - M^N_t\rangle\,dt\right] + \frac{C'}{c}\sup_{\varphi \in \mathcal{T}_c}d_2(\varphi, \tilde{\mathcal{T}}_c) + (8T\delta_N)^{1/2}\right)\\
            &\leq C\left(\frac{1}{c}\sup_{\tilde{\varphi} \in \tilde{\mathcal{T}}_c}L_N(\lambda^N, \tilde{\varphi}) + \frac{C'}{c}\sup_{\varphi \in \mathcal{T}_c}d_2(\varphi, \tilde{\mathcal{T}}_c) + (8T\delta_N)^{1/2}\right)
        \end{align*}
        Thus by \cref{lemma:d_2_univ_approx_error} and after pulling out constants and relabeling we have that:
        \begin{align*}
            d_2(\lambda^N, \lambda^*) \leq C\left(\sup_{\tilde{\varphi} \in \tilde{\mathcal{T}}_c}L_N(\lambda^N, \tilde{\varphi}) + \sup_{\varphi \in \mathcal{D}}d_2(\varphi, \tilde{\mathcal{D}}) + \delta_N^{1/2}\right), \quad \forall \lambda^N \in \mathcal{A}, \ ||\lambda^N||_{\mathbb{S}^2_{\mathbb{F}}} \leq C_{\lambda}
        \end{align*}
        Now we have that:
        \begin{align*}
            L_N(\lambda^N, \tilde{\varphi}) &= \mathbb{E}\left[||g(v^N_{t_N}) - \lambda^N_{t_N}||_2^2\right]^{1/2} + \mathbb{E}\left[\int_0^T\langle \tilde{\varphi}_t, Y^N_t - M^N_t\rangle\,dt\right] \\
            &= \mathbb{E}\left[||g(v^N_{t_N}) + \int_0^Tb(\tau(s), v_{\tau(s)}^N)\,ds -  M^N_{t_N}||_2^2\right]^{1/2} + \mathbb{E}\left[\int_0^T\langle \tilde{\varphi}_t, \tilde{M}^N_t - M^N_t\rangle\,dt\right]\\
            &\leq \mathbb{E}\left[||\tilde{M}^N_{t_N} -  M^N_{t_N}||_2^2\right]^{1/2} + ||\tilde{\varphi}||_{\mathbb{L}^2_{\mathbb{F}}}||\tilde{M}^N - M^N||_{\mathbb{L}^2_{\mathbb{F}}} \\
            &\leq (1 + c\sqrt{T})d_S(\tilde{M}^N, M^N) 
        \end{align*}
        Where recalling \cref{remark:tilde_m} we have $\tilde{M}^N \in \mathbb{M}^{2}_{\mathbb{F}}([0,T]; \mathbb{R}^d)$. Thus we have:
        \begin{align*}
            \inf_{\lambda^N \in \mathscr{A}_{\lambda}}\sup_{\tilde{\varphi} \in \tilde{\mathcal{T}}_c}L_N(\lambda^N, \tilde{\varphi}) \leq (1+c\sqrt{T})\sup_{M \in \mathbb{M}^{2,*}_{\mathbb{F}}}\inf_{\lambda^N \in \mathscr{A}_{\lambda}}d_S(M, M^N(\lambda^N))
        \end{align*}
        proving \ref{eq:discrete_min_max_lag_bound}.
    \end{proof}
    \begin{proof}[Proof of 2]
        Letting:
        \begin{align*}
            X_t = \mathbb{E}\left[g(v_T) - g(v_{t_N}^N) + \int_t^T(b(s, v_s) - b(\tau(s),v_{\tau(s)}^N)\,ds \ \bigg| \mathcal{F}_t\right], \quad t \in [0,T]
        \end{align*}
        for $\varphi \in  \mathcal{D}$ by \cref{lemma:discrete_fpr} we have that:
        \begin{align*}
            \mathbb{E}\left[\int_0^T\langle \varphi_t, \Phi(\lambda^N)_t - \lambda^N_t\rangle\,dt\right] = \mathbb{E}\left[\int_0^T\langle \varphi_t, \tilde{M}_t^N - M^N_t\rangle\,dt\right] 
            + \mathbb{E}\left[\int_0^T\langle \varphi_t, X_t \rangle\,dt\right] 
        \end{align*}
        thus we have:
        \begin{align*}
             &\mathbb{E}\left[\int_0^T\langle \varphi_t, \Phi(\lambda^N)_t - \lambda^N_t\rangle\,dt\right] -\frac{r}{2}||\varphi||^2_{\mathbb{L}^2_{\mathbb{F}}} = \mathbb{E}\left[\int_0^T\langle \tilde{\varphi}_t, Y_t^N - M^N_t\rangle\,dt\right] -  \frac{r}{2}||\tilde{\varphi}||^2_{\mathbb{L}^2_{\mathbb{F}}} \\
             &\quad\quad +\mathbb{E}\left[\int_0^T\langle \varphi_t - \tilde{\varphi}_t, Y_t^N - M^N_t\rangle\,dt\right] + \frac{r}{2}\left( ||\tilde{\varphi}||^2_{\mathbb{L}^2_{\mathbb{F}}}-||\varphi||^2_{\mathbb{L}^2_{\mathbb{F}}} \right) \\
             &\quad\quad + \mathbb{E}\left[\int_0^T\langle \varphi_t, X_t \rangle\,dt\right] 
        \end{align*}
        By Cauchy-Schwarz, the triangle inequality and \cref{lemma:bounded_Y_M_min_max_solver}, and 
        using $||a||^2 - ||b||^2 = \langle a - b, a+b\rangle \leq ||a-b||||a+b||$ we have for $(\varphi, \tilde{\varphi}) \in \mathcal{D}\times \tilde{\mathcal{D}}$ that:
        \begin{align*}
            &\mathbb{E}\left[\int_0^T\langle \varphi_t, \Phi(\lambda^N)_t - \lambda^N_t\rangle\,dt\right] - \frac{r}{2}||\varphi||^2_{\mathbb{L}^2_{\mathbb{F}}} \leq \mathbb{E}\left[\int_0^T\langle \tilde{\varphi}_t, Y_t^N - M^N_t\rangle\,dt\right] -  \frac{r}{2}||\tilde{\varphi}||^2_{\mathbb{L}^2_{\mathbb{F}}} \\
            &\quad\quad + ||\varphi - \tilde{\varphi}||_{\mathbb{L}^2_{\mathbb{F}}}\left(C'+ \frac{r}{2}||\tilde{\varphi} + \varphi||_{\mathbb{L}^2_{\mathbb{F}}}\right) + ||\varphi||_{\mathbb{L}^2_{\mathbb{F}}}\mathbb{E}\left[\int_0^T||X_t||_2^2\,dt\right]^{1/2}
        \end{align*}
        Now defining $\delta_N \geq 0$ such that:
        \begin{align*}
            \delta_N := \sup_{\lambda^N \in \mathcal{A}}\mathbb{E}\left[||g(v_T) - g(v_{t_N}^N)||^2_2 + T\int_0^T||b(s, v_s) - b(\tau(s),v_{\tau(s)}^N)||_2^2\,ds\right]
        \end{align*}
        we have by the definition of $\tilde{L}_N$, the triangle inequality and \cref{lemma:discrete_fpr} that:
        \begin{align*}
            &\mathbb{E}\left[\int_0^T\langle \varphi_t, \Phi(\lambda^N)_t - \lambda^N_t\rangle\,dt\right] - \frac{r}{2}||\varphi||^2_{\mathbb{L}^2_{\mathbb{F}}} \leq \sup_{\tilde{\varphi} \in \tilde{\mathcal{D}}}\tilde{L}_N(\lambda^N, \tilde{\varphi}) \\
            &\quad\quad + ||\varphi - \tilde{\varphi}||_{\mathbb{L}^2_{\mathbb{F}}}\left(C'+ \frac{r}{2}||\tilde{\varphi}||_{\mathbb{L}^2_{\mathbb{F}}} + \frac{r}{2}||\varphi||_{\mathbb{L}^2_{\mathbb{F}}}\right)  + ||\varphi||_{\mathbb{L}^2_{\mathbb{F}}}(8T\delta_N)^{1/2}, \quad \forall (\varphi, \tilde{\varphi}) \in \mathcal{D}\times \tilde{\mathcal{D}}
        \end{align*}
        Now for $c > 0$ defining the sets of bounded test functions:
        \begin{align*}
             \mathcal{T}_c:= \{\varphi \in \mathcal{D} : ||\varphi||_{\mathbb{L}^2_{\mathbb{F}}} \leq c\}, \quad \tilde{\mathcal{T}}_c := \tilde{\mathcal{D}}\cap \mathcal{T}_c
        \end{align*}
        We have that:
        \begin{align*}
            &\mathbb{E}\left[\int_0^T\langle \varphi_t, \Phi(\lambda^N)_t - \lambda^N_t\rangle\,dt\right] - \frac{r}{2}||\varphi||^2_{\mathbb{L}^2_{\mathbb{F}}} \leq \sup_{\tilde{\varphi} \in \tilde{\mathcal{D}}}\tilde{L}_N(\lambda^N, \tilde{\varphi}) \\
            &\quad\quad + ||\varphi - \tilde{\varphi}||_{\mathbb{L}^2_{\mathbb{F}}}\left(C'+rc\right) + c(8T\delta_N)^{1/2}, \quad \forall (\varphi, \tilde{\varphi}) \in \mathcal{T}_c\times \tilde{\mathcal{T}}_c
        \end{align*}
        Now defining:
        \begin{align*}
            d_2(\varphi, \tilde{\mathcal{T}}_c) := \inf_{\tilde{\varphi} \in \tilde{\mathcal{T}}_c}d_2(\varphi, \tilde{\varphi}), \quad \varphi \in \mathcal{T}_c
        \end{align*}
        where by \cref{lemma:d_2_univ_approx_error}, we have that:
        \begin{align*}
            \sup_{\varphi \in \mathcal{T}_c}d_2(\varphi, \tilde{\mathcal{T}}_c) \leq \sup_{\varphi \in \mathcal{D}}d_2(\varphi, \tilde{\mathcal{D}}) 
        \end{align*}
        thus we obtain:
        \begin{align*}
            &\mathbb{E}\left[\int_0^T\langle \varphi_t, \Phi(\lambda^N)_t - \lambda^N_t\rangle\,dt\right] - \frac{r}{2}||\varphi||^2_{\mathbb{L}^2_{\mathbb{F}}} \leq \sup_{\tilde{\varphi} \in \tilde{\mathcal{D}}}\tilde{L}_N(\lambda^N, \tilde{\varphi}) \\
            &\quad\quad + \left(C'+rc\right)\sup_{\varphi \in \mathcal{D}}d_2(\varphi, \tilde{\mathcal{D}}) + c(8T\delta_N)^{1/2}, \quad \forall \varphi \in \mathcal{T}_c
        \end{align*}
        Now by standard arguments we have that:
        \begin{align*}
            \frac{1}{2r}d_2(\Phi(\lambda^N), \lambda^N)^2 &= \sup_{\varphi \in \mathcal{D}}\left\{ \mathbb{E}\left[\int_0^T\langle \varphi_t, \Phi(\lambda^N)_t - \lambda^N_t\rangle\,dt\right] - \frac{r}{2}||\varphi||^2_{\mathbb{L}^2_{\mathbb{F}}}\right\} \\ 
            &= \mathbb{E}\left[\int_0^T\langle \varphi_t^*, \Phi(\lambda^N)_t - \lambda^N_t\rangle\,dt\right] - \frac{r}{2}||\varphi^*||^2_{\mathbb{L}^2_{\mathbb{F}}} 
        \end{align*}
        Where we have $\varphi^* = (\Phi(\lambda^N) - \lambda^N)/ r \in \mathbb{L}^2_{\mathbb{F}}([0,T]; \mathbb{R}^d)$, possibly not in $\mathcal{D}$. Since $||\lambda^N||_{\mathbb{L}^2_{\mathbb{F}}} \leq \sqrt{T}C_{\lambda}$, by \cref{lemma:bounded_Y_M_min_max_solver} there exists $C'' > 0$ such that $||\varphi^*||_{\mathbb{L}^2_{\mathbb{F}}} \leq (C'' + \sqrt{T}C_{\lambda})/r$. Choosing $c > (C'' + \sqrt{T}C_{\lambda})/r$ we have that:
        \begin{align*}
            \sup_{\varphi \in \mathcal{D}}\left\{ \mathbb{E}\left[\int_0^T\langle \varphi_t, \Phi(\lambda^N)_t - \lambda^N_t\rangle\,dt\right] - \frac{r}{2}||\varphi||^2_{\mathbb{L}^2_{\mathbb{F}}}\right\} = \sup_{\varphi \in \mathcal{T}_c}\left\{ \mathbb{E}\left[\int_0^T\langle \varphi_t, \Phi(\lambda^N)_t - \lambda^N_t\rangle\,dt\right] - \frac{r}{2}||\varphi||^2_{\mathbb{L}^2_{\mathbb{F}}}\right\}
        \end{align*}
        Consequently we have that:
        \begin{align*}
             \frac{1}{2r}d_2(\Phi(\lambda^N),\lambda^N) \leq \sup_{\tilde{\varphi} \in \tilde{\mathcal{D}}}\tilde{L}_N(\lambda^N, \tilde{\varphi}) + \left(C'+rc\right)\sup_{\varphi \in \mathcal{D}}d_2(\varphi, \tilde{\mathcal{D}}) + c(8T\delta_N)^{1/2}
        \end{align*}
        Now by \cref{lemma:fixed_point_stability} there exists $C > 0$ such that $d_2(\lambda, \lambda^*)^2 \leq Cd_2(\Phi(\lambda), \lambda)^2, \forall \lambda \in \mathcal{A}$, so we obtain:
        \begin{align*}
            d_2(\lambda^N), \lambda^*)^2  \leq 2rC\left(\sup_{\tilde{\varphi} \in \tilde{\mathcal{D}}}\tilde{L}_N(\lambda^N, \tilde{\varphi}) + \left(C'+rc\right)\sup_{\varphi \in \mathcal{D}}d_2(\varphi, \tilde{\mathcal{D}}) + c(8T\delta_N)^{1/2}\right)
        \end{align*}
        Thus after pulling out constants and relabeling we obtain:
        \begin{align*}
            d_2(\lambda^N, \lambda^*)^2  \leq C\left(\sup_{\tilde{\varphi} \in \tilde{\mathcal{D}}}\tilde{L}_N(\lambda^N, \tilde{\varphi})
            +\sup_{\varphi \in \mathcal{D}}d_2(\varphi, \tilde{\mathcal{D}})  +\delta_N^{1/2}\right)
        \end{align*}      
        Now, we have by the tower property and standard arguments that:
        \begin{align*}
            \tilde{L}_N(\lambda^N, \tilde{\varphi}) &= \mathbb{E}\left[||g(v^N_{t_N}) - \lambda^N_{t_N}||_2^2\right] + \mathbb{E}\left[\int_0^T\langle \tilde{\varphi}_t, \tilde{M}^N_t - M^N_t\rangle\,dt\right] - \frac{r}{2}||\tilde{\varphi}||_{\mathbb{L}^2_{\mathbb{F}}}^2\\
            &\leq \mathbb{E}\left[||g(v^N_{t_N}) - \lambda^N_{t_N}||_2^2\right] + \frac{1}{2r}\mathbb{E}\left[\int_0^T|| \tilde{M}^N_t - M^N_t||_2^2\,dt\right]  \\
            &= \mathbb{E}\left[||g(v^N_{t_N}) + \int_0^Tb(\tau(s), v_{\tau(s)}^N)\,ds -  M^N_{t_N}||_2^2\right] + \frac{1}{2r}\mathbb{E}\left[\int_0^T|| \tilde{M}^N_t - M^N_t||_2^2\,dt\right] \\
            &\leq \mathbb{E}\left[||\tilde{M}^N_{t_N} -  M^N_{t_N}||_2^2\right] + \frac{1}{2r}||\tilde{M}^N - M^N||_{\mathbb{L}^2_{\mathbb{F}}}^2 \\
            &\leq \left(1 + \frac{T}{2r}\right)d_S(\tilde{M}^N, M^N)^2
        \end{align*}
        Where we have $\tilde{M}^N \in \mathbb{M}^{2}_{\mathbb{F}}([0,T]; \mathbb{R}^d)$. Thus we obtain:
        \begin{align*}
            \inf_{\lambda^N \in \mathscr{A}_{\lambda}}\sup_{\tilde{\varphi} \in \tilde{\mathcal{D}}}\tilde{L}_N(\lambda^N, \tilde{\varphi}) \leq \left(1+\frac{T}{2r}\right)\sup_{M \in \mathbb{M}^{2}_{\mathbb{F}}}\inf_{\lambda^N \in \mathscr{A}_{\lambda}}d_S(M, M^N(\lambda^N))^2
        \end{align*}
        proving \ref{eq:discrete_min_max_reg_lag_bound}.
    \end{proof}
    
\section{Conclusion}
In this paper we have developed two solvers (a martingale solver, and a Hilbert space orthogonal min-max solver) for the coupled forward-backward system of \cref{eq:fbs_general} using the useful equivalent fixed point formulation developed in \cite{KeyworthPaper:2026}. Proving an important fixed point stability result for this problem enabled us to prove a posteriori error bounds for both solvers. We studied the properties of discretized versions of both solvers and again used the fixed point stability result to obtain a posteriori error bounds subject to error terms coming from approximation error (controlled by the choice of approximating class of functions), and time discretization error. The numerical implementation of the martingale based solver would require the class of square integrable martingales to be parameterized in some way. For deep FBSDE solvers previously studied in the literature where the backward equation is a BSDE with a Brownian driver, the parameterization is achieved by approximating the initial condition and integrand of the stochastic integral using neural networks. This method can be applied in the case of our coupled forward-backward system if the filtration possess the martingale representation property with respect to a collection of known martingales that can be sampled. Our Hilbert space orthogonal solver does not require any such restriction and is suitable for general filtrations. The key distinction between the two solvers lies in the topology in which the solver objective controls the solution error. In particular, given two martingales, the martingale property ensures that controlling only the difference between their terminal values suffices to control their difference over the entire path, uniformly in time. The martingale solver inherits this property and is therefore stable in a strong topology that controls whole path differences uniformly in time. However, the min-max solvers provide stability only in a weaker $L^2$-type topology, which follows naturally from the fact that we rely on a Hilbert space orthogonality condition to try and enforce the martingale property.

\bibliographystyle{alpha}
\bibliography{references}

\appendix

\section{Auxiliary Results}
\begin{lemma}\label{lemma:dense_test_l2}
    Let $\mathcal{H}$ be a Hilbert space, and $\Phi:\mathcal{H}\mapsto \mathbb{R}$ a continuous homogeneous functional. Let $\mathcal{D}\subseteq\mathcal{H}$ be a dense linear subspace. Given $c > 0$, if $\Phi(h) = 0, \ \forall h \in \{h \in \mathcal{D}: ||h||\leq c\}$ then $\Phi(h) = 0, \forall h \in \mathcal{H}$.
    \begin{proof}
        By assumption $\Phi(0) = 0$. For any $h \in \mathcal{D}$ with $||h||\neq 0$ we have that $g := ch/||h|| \in \{h \in \mathcal{D}: ||h||\leq c\}$, thus by homogeneity $\Phi(h) = \Phi(g)||h||/c = 0$ and hence $\Phi(h)=0, \forall h \in \mathcal{D}$. Since $\mathcal{D}\subseteq \mathcal{H}$ is dense, for any $h \in \mathcal{H}$ there exists a sequence $(h_n)_{n \in \mathbb{N}}\subseteq \mathcal{D}$ with $h_n\rightarrow h$. By continuity we have that:
        \begin{align*}
            \lim_{n\rightarrow \infty}\Phi(h_n) = \Phi(h)
        \end{align*}
        but $\Phi(h_n) = 0, \forall n \in \mathbb{N}$, thus $\lim_{n\rightarrow\infty}\Phi(h_n) = \Phi(h) = 0$. Since $h$ was arbitrary the result follows.
    \end{proof}
\end{lemma}
\begin{lemma} \label{lemma:d_2_univ_approx_error}
    Let $(\mathcal{H}, \langle \cdot, \cdot\rangle)$ be a Hilbert space and let $\mathcal{D}'\subseteq \mathcal{D}\subseteq \mathcal{H}$ with $\mathcal{D}$ a dense subset and $\mathcal{D}'$ a linear subspace. For $c > 0$ define the bounded set $\mathcal{T}_c := \{x \in \mathcal{D} : ||x||\leq c\}$ and $\mathcal{T}':=\mathcal{D}'\cap \mathcal{T}$. Then:
    \begin{align*}
        \sup_{x \in \mathcal{T}}d(x, \mathcal{T}') \leq 2\sup_{x \in \mathcal{D}}d(x,\mathcal{D}')
    \end{align*}
    where for $A\subseteq \mathcal{H}$ and $x \in \mathcal{H}$ we define $d(x,A) := \inf_{y \in A}||x-y||$.
    \begin{proof}
        We have that:
        \begin{align*}
            d(x, \mathcal{T}') = \inf_{x' \in \mathcal{T}'}||x - x'|| \leq ||x - x'||, \quad \forall x, x' \in \mathcal{D}\times \mathcal{T}'
        \end{align*}
        For any $x' \in \mathcal{D}'$ we have $\tilde{x} := cx'/||x'|| \in \mathcal{T}'$ since $\mathcal{D}'$ is a linear subspace. Thus:
        \begin{align*}
            d(x, \mathcal{T}') \leq ||x - \tilde{x}|| \leq ||x - x'||+||x' - \tilde{x}|| \leq 2||x - x'||, \quad \forall (x, x') \in \mathcal{T}\times \mathcal{D}'
        \end{align*}
        by the fact that:
        \begin{align*}
            ||x' - \tilde{x}|| = ||x'||(1 - c/||x'||) = ||x'|| - c \leq ||x'|| - ||x|| \leq ||x - x'||
        \end{align*}
        Thus we have:
        \begin{align*}
            \sup_{x \in \mathcal{T}}d(x, \mathcal{T}') \leq 2\sup_{x \in \mathcal{T}}d(x, \mathcal{D}') \leq 2\sup_{x \in \mathcal{D}}d(x, \mathcal{D}')
        \end{align*}
    \end{proof}
\end{lemma}

\begin{lemma} \label{lemma:G_adaptedness}
    Let $i$ be the canonical isometric embedding such that:
    \begin{align*}
        \mathbb{L}^2_{\mathbb{F}}([0,T]; \mathbb{R}^d)  \xhookrightarrow{i} \mathbb{L}^2([0,T]; \mathbb{R}^d)
    \end{align*}
    defined by:
    \begin{align*}
        i : \mathbb{L}^2_{\mathbb{F}}([0,T]; \mathbb{R}^d)\rightarrow \mathbb{L}^2([0,T]; \mathbb{R}^d), \quad i([X]_{d_2}) = [X]_{d_2}
    \end{align*}
    We have:
    \[\forall [X]_{d_2} \in i(\mathbb{L}^2_{\mathbb{F}}([0,T]; \mathbb{R}^d)), \quad \forall Y \in [X]_{d_2}, \quad Y(t,\cdot) \text{ is } \mathcal{F}_t\text{-measurable a.e. $t \in [0,T]$} \]
    \begin{proof}
        For $[X]_{d_2} \in i(\mathbb{L}^2_{\mathbb{F}})$ there exists a $\mathbb{F}$-adapted representative $X \in [X]_{d_2}$. By definition and Fubini: $X_t= Y_t, \ \mathbb{P}\text{-a.s.}, \text{ a.e. } t \in [0,T], \ \forall Y \in [X]_{d_2}$.
        By completeness of $\mathcal{F}_t$ it follows that $Y_t$ is $\mathcal{F}_t$ measurable a.e. $t \in [0,T]$. 
    \end{proof}
\end{lemma}
\begin{lemma} \label{lemma:orth_proj}
    Let $i$ be the canonical isometric embedding such that:
    \begin{align*}
        \mathbb{L}^2_{\mathbb{F}}([0,T]; \mathbb{R}^d)  \xhookrightarrow{iso} \mathbb{L}^2([0,T]; \mathbb{R}^d)
    \end{align*}
    Let $\pi:\mathbb{L}^2\rightarrow i(\mathbb{L}^2_{\mathbb{F}})$ denote the orthogonal projection operator, then we have that:
    \begin{equation}
        \begin{aligned}
            (\pi X)_t &= \mathbb{E}\left[X_t|\mathcal{F}_t\right], \quad \mathbb{P}\text{-a.s.}, \ \text{a.e. } t \in [0,T], \quad \forall X \in \mathbb{L}^2 \\
            (\pi X)_t &= X_t, \quad \mathbb{P}\text{-a.s.}, \ \text{a.e. } t \in [0,T], \quad \forall X \in i(\mathbb{L}^2_{\mathbb{F}})
        \end{aligned}
    \end{equation}
    \begin{proof}
        For $X \in \mathbb{L}^2$ we have $X - \pi X \in i(\mathbb{L}^2_{\mathbb{F}})^\perp$ and thus:
        \begin{align*}
            \langle h, X - \pi X\rangle_{\mathbb{L}^2} = \mathbb{E}\left[\int_0^T\langle h_t,  X_t -(\pi X)_t\rangle\,dt\right] = 0, \quad \forall h \in i(\mathbb{L}^2_{\mathbb{F}})
        \end{align*}
        Thus by Fubini we have that:
        \begin{align*}
            \mathbb{E}\left[\langle h_t,  X_t -(\pi X)_t\rangle\right] = 0, \quad \forall h \in i(\mathbb{L}^2_{\mathbb{F}}), \quad \text{a.e. } t \in [0,T] 
        \end{align*}
        By \cref{lemma:G_adaptedness} we have that $h_t, (\pi X)_t \in L^2(\Omega, \mathcal{F}_t, \mathbb{P}; \mathbb{R}^d), \ \text{a.e. } t \in [0,T]$ thus by the orthogonal projection theorem on $L^2(\Omega, \mathcal{F}_t, \mathbb{P}; \mathbb{R}^d)$ we have that:
        \begin{align*}
            (\pi X)_t &= \mathbb{E}\left[X_t|\mathcal{F}_t\right], \quad \mathbb{P}\text{-a.s.}, \ \text{a.e. } t \in [0,T]
        \end{align*}
        Now supposing that $X \in i(\mathbb{L}^2_{\mathbb{F}})$, \cref{lemma:G_adaptedness} we have that:
        \begin{align*}
            (\pi X)_t = \mathbb{E}[X_t|\mathcal{F}_t] = X_t, \quad \mathbb{P}\text{-a.s.}, \ \text{a.e. } t \in [0,T]
        \end{align*}
    \end{proof}
\end{lemma}

\begin{lemma} \label{corr:orth_proj_cost}
    Let $i$ be the canonical embedding such that:
    \begin{align*}
        \mathbb{L}^2_{\mathbb{F}}([0,T]; \mathbb{R}^d)  \xhookrightarrow{iso} \mathbb{L}^2([0,T]; \mathbb{R}^d)
    \end{align*}
    and let $\mathcal{D}\subseteq i(\mathbb{L}^2_{\mathbb{F}})$ be a dense linear subspace and for any $c > 0$ define the set $\mathcal{T}_c:=\{\varphi \in \mathcal{D}: ||\varphi|| \leq c\}$. Let $X \in \mathcal{A}$ and  define $Y\in \mathbb{L}^2$ such that $(Y_t)_{t \in [0,T]} \equiv Z \in L^2(\Omega, \mathcal{F}_T, \mathbb{P}; \mathbb{R}^d)$. Then:
    \begin{align*}
        \sup_{\varphi \in \mathcal{T}_c}\langle g, Y - X\rangle_{\mathbb{L}^2} = 0
    \end{align*}
    if and only if $X$ is a version of the conditional expectation process:
    \begin{align*}
        X_t = \mathbb{E}[Z|\mathcal{F}_t], \quad \mathbb{P}\text{-a.s.}, \quad t \in [0,T]
    \end{align*}
    \begin{proof}
       Define the continuous homogeneous functional $F: i(\mathbb{L}^2_{\mathbb{F}})\rightarrow \mathbb{R}$:
        \begin{align*}
            F(\varphi) := \langle \varphi, Y - X\rangle_{\mathbb{L}^2}
        \end{align*}
        Suppose that:
        \begin{align*}
            \sup_{\varphi \in \mathcal{T}_c}F(\varphi) = 0, \quad \exists \varphi \in \mathcal{T}_c : F(\varphi) < 0
        \end{align*}
        Then we have $-\varphi \in \mathcal{T}_c$ and $F(-\varphi) > 0$ which is a contradiction. Thus:
        \begin{align*}
            \sup_{\varphi \in \mathcal{T}_c}F(\varphi) = 0 \implies F(\varphi) = 0, \ \forall \varphi \in \mathcal{T}_c
        \end{align*}
        By \cref{lemma:dense_test_l2} if $F(\varphi) = 0, \forall \varphi \in \mathcal{T}_c$ then $F(\varphi) = 0, \forall \varphi \in i(\mathbb{L}^2_{\mathbb{F}})$. Thus we have that $Y - X \in i(\mathbb{L}^2_{\mathbb{F}})^\perp$, and thus:
        \begin{align*}
            \pi (Y - X) = \pi Y - \pi X = 0 \in i(\mathbb{L}^2_{\mathbb{F}})
        \end{align*}
        and by \cref{lemma:orth_proj} we obtain:
        \begin{align*}
            \sup_{\varphi \in \mathcal{T}_c}F(\varphi) = 0 \iff X_t = \mathbb{E}[Z|\mathcal{F}_t], \quad \mathbb{P}\text{-a.s.}, \quad a.e. \ t \in [0,T]
        \end{align*}
        The square integrable martingale $(\mathbb{E}[Z|\mathcal{F}_t])_{t\in [0,T]}$ on the RHS admits a right continuous version $(M_t)_{t \in [0,T]}$ (which exists by the corollary of Theorem 3.1 in \cite{LipsterBook1:2001}). Hence we must have that:
        \begin{align*}
            X_t = M_t, \quad \mathbb{P}\text{-a.s.}, \quad \text{a.e.} \ t \in [0,T]
        \end{align*}
        But $X$ is right continuous, and two right continuous processes $\mathbb{P}\text{-a.s.}$ equal $\text{a.e. } t \in [0,T]$ are indistinguishable by standard arguments. 
    \end{proof}
\end{lemma}

\end{document}